\documentclass[a4paper,12pt]{article}

\newcommand{\plim}[1][]{\mathop{\varprojlim}\limits_{#1}}
\newcommand{\ilim}[1][]{\mathop{\varinjlim}\limits_{#1}}
\usepackage{amssymb}
\usepackage{mathrsfs}
\usepackage{amsmath}
\usepackage{amsfonts}
\allowdisplaybreaks[4]
\usepackage{amsthm}
\usepackage{amscd}
\usepackage[all]{xy}
\usepackage{comment}
\usepackage{fancybox}

\makeatletter
    
    \@addtoreset{equation}{section}
  \makeatother

\newtheoremstyle{break}
  {}{}{\itshape}{}{\bfseries}{.}{\newline}{}

\theoremstyle{break}
\newtheorem{thm}{Theorem}[section]
\newtheorem{lem}{Lemma}[section]
\newtheorem{cor}{Corollary}[section]
\newtheorem{prop}{Proposition}[section]
\newtheorem*{prop*}{Proposition}
\newtheorem{ex}{Example}[section]

\newtheoremstyle{nonbreak}
  {}{}{\itshape}{}{\bfseries}{.}{ }{}
\theoremstyle{nonbreak}
\newtheorem{nonthm}[thm]{Theorem}
\newtheorem{nonprop}[prop]{Proposition}

\theoremstyle{definition}
\newtheorem{defn}{Definition}[section]
\newtheorem*{rem}{Remark}
\newtheorem*{prf}{Proof}

\newtheorem*{cnj*}{Conjecture}

\DeclareMathOperator{\rank}{rank}
\DeclareMathOperator{\Ker}{Ker}
\DeclareMathOperator{\Image}{Im}
\DeclareMathOperator{\Real}{Re}
\DeclareMathOperator{\Hom}{Hom}
\DeclareMathOperator{\Ext}{Ext}
\DeclareMathOperator{\supp}{supp}
\DeclareMathOperator{\Tor}{Tor}

\title{Residue current approach to Ehrenpreis-Malgrange type theorem for linear differential equations with constant coefficients and commensurate time lags}

\author{Saiei-Jaeyeong Matsubara-Heo\footnote{saiei@ms.u-tokyo.ac.jp}}

\begin{document}

\date{}
\maketitle

\begin{abstract}
 We introduce a ring $\mathcal{H}$ of partial difference-differential operators with constant coefficients initially defined by H. Gl\"using-L\"ur\ss en for ordinary difference-differential operators and investigate its cohomological properties.
Combining this ring theoretic observation with the integral representation technique developed by M. Andersson, we solve a certain type of division with bounds. In the last section, we deduce from this injectivity properties of various function modules over $\mathcal{H}$ as well as the density results of exponential polynomial solutions for partial difference-differential equations.
\end{abstract}
%ここは全数学者に理解できるように書くべき。

%\begin{keyword}
%Residue currents, Ehrenpreis' fundamental principle, Difference-Differential equations, Division with bounds, Coherent rings


%\end{keyword}

%\end{frontmatter}



\section{Introduction}

%\fbox{Emphasize that $\mathcal{H}$ is the best possible ring!!}
In the last century, various studies on the general theory of linear partial differential equations with constant coefficients has been conducted. Among many important works, we should mention the celebrated work of L. Ehrenpreis: Ehrenpreis' fundamental principle. Roughly speaking, Ehrenpreis' fundamental principle means that any solution of a system of linear partial differential equations with constant coefficients can be written as a superposition of exponential polynomial solutions. In particular, exponential polynomial solutions are dense in the solution space.\footnote{This type of theorem can never be true in variable coefficients cases. In this paper,  ``difference-differential equations'' always means linear difference-differential equations with constant coefficients.} Furthermore, if we let $\mathbb{C}[\partial]=\mathbb{C}[\partial_1,\cdots,\partial_n]$ denote the ring of linear partial differential equations with constant coefficients, Ehrenpreis fundamental principle has its cohomological counterpart which is expressed as injectivity properties of function spaces over $\mathbb{C}[\partial]$. We give a statement of a fundamental result in the theory of linear partial differential equations only for $C^\infty$ functions.

\begin{thm}[{\cite[Theorem 4.2]{Eh}}, {\cite[Theorem 7.6.12, 7.6.13]{Ho}}, {\cite[IV,\S 5,$5^\circ$]{P}}]\label{EM}
Let $\Omega$ be a convex subset of $\mathbb{R}^n$ and let $\mathbf{E}$ denote the $\mathbb{C}[\partial]$-module consisting of exponential polynomials. One has the following properties:
\begin{enumerate}
\item $C^\infty(\Omega)$ is an injective $\mathbb{C}[\partial]$-module, i.e., for any finitely generated $\mathbb{C}[\partial]$-module $M$ and any positive integer $i,$ we have an identity
\begin{equation}
\Ext_{\mathbb{C}[\partial]}^i(M,C^\infty(\Omega))=0.
\end{equation}
\item For any matrix $\mathbb{P}(\partial)\in M(r_1,r_0;\mathbb{C}[\partial])$, 
$\Ker\left(\mathbb{P}(\partial):\mathbf{E}^{r_0\times 1}\rightarrow\mathbf{E}^{r_1\times 1}\right)$ is dense in $\Ker(\mathbb{P}(\partial):C^\infty(\Omega)^{r_0\times 1}\rightarrow C^\infty(\Omega)^{r_1\times 1}).$
\end{enumerate}
\end{thm}

We call the theorem above Ehrenpreis-Malgrange type theorem. In the following, we explain the meaning of (1) of Theorem \ref{EM}. To this end, let us consider a more general setting.

%Our main result is a generalization of this density result to difference-differential equations:

Let $R$ be a subring of the ring of distributions with compact supports $\mathcal{E}^\prime(\mathbb{R}^n)$ and $F$ be a vector space of functions of a suitable kind, say $F=C^\infty(\mathbb{R}^n)$ for simplicity. Then, $R$ acts on $F$ by 
\begin{equation}
T\cdot f\overset{def}{=}T*f\hspace{3mm}(T\in R,f\in F).
\end{equation}
Here, $*$ is the usual convolution product. Note that usual partial differential operator $\frac{\partial}{\partial x_i}$ %and a difference operator $\sigma$ to the first coordinate
 can be realized as a convolution operator by 
\begin{equation}
\frac{\partial f}{\partial x_i}=\left(\frac{\partial}{\partial x_i}\delta\right) *f,
\end{equation}
%\text{ ~resp.~ }\sigma f=f(x_1+1,x_2,\cdots,x_n)=\delta(x_1+1,x_2,\cdots,x_n)*f,$$
where $\delta$ is Dirac measure with support at the origin. For any $r_1\times r_0$ matrix $\mathbb{P}$ with entries in $R$ and for any $\mathbf{f}={}^t(f_1,\cdots,f_{r_1})\in F^{r_1\times 1}$, we consider a system of equations 
\begin{equation}\label{TheEquation0}
\mathbb{P}\mathbf{u}=\mathbf{f},
\end{equation}
where $\mathbf{u}\in F^{r_0\times 1}$ is an unknown vector of functions. Now, suppose this system of equations (\ref{TheEquation0}) has a solution $\mathbf{u}.$ Then, we can easily find constraints on $\mathbf{f}$:
\begin{equation}\label{Compatibility}
\text{For any }r_2\times r_1\text{ matrix }\mathbb{Q}\text{ with entries in } R, \;\mathbb{Q}\mathbb{P}=0\text{ implies } \mathbb{Q}\mathbf{f}=0.
\end{equation} The condition (\ref{Compatibility}) is called the compatibility condition of the system (\ref{TheEquation0}). One may ask if this condition is sufficient for the solvability of the system (\ref{TheEquation0}). After a simple  discussion on homological algebra, the problem is reduced to proving the identity

\begin{equation}
\Ext_R^1(R^{1\times r_0}/R^{1\times r_1}\mathbb{P},F)=0.
\end{equation}
Thus, by (1) of Theorem \ref{EM}, the system (\ref{TheEquation0}) is solvable if and only if (\ref{Compatibility}) is satisfied when $R=\mathbb{C}[\partial]$.

Since the formulation above is quite general, it is natural to ask if one can extend  Ehrenpreis-Malgrange type theorem to more general convolution equations, i.e., if one can replace $\mathbb{C}[\partial]$ by a larger ring $R$ in Theorem \ref{EM}. For example, consider a system

\begin{equation}
\mathbb{P}(\partial,\sigma)\mathbf{u}=\mathbf{f},
\end{equation}
where $\mathbb{P}(\partial,\sigma)$ is a $r_1\times r_0$ matrix with entries in $\mathbb{C}[\partial,\sigma]=\mathbb{C}[\partial_1,\cdots,\partial_n,\sigma].$ Here, $\sigma$ is a difference operator to the first coordinate which can be realized as a convolution operator by
\begin{equation}
\sigma f=f(x_1+1,x_2,\cdots,x_n)=\delta(x_1+1,x_2,\cdots,x_n)*f.
\end{equation}
When we regard the first coordinate as a time variable, this system of difference-differential equations is called linear partial differential equations with constant coefficients and commensurate time lags. For brevity, we call such a system $D\Delta$-equations in this paper. This system has attracted attentions of many mathematicians, especially in connection with ring theory (\cite{OB}, \cite{GL}). One might hope that Ehrenpreis-Malgrange type theorem holds for $R=\mathbb{C}[\partial,\sigma].$ However, this optimistic conjecture turns out to be false even if $n=1$ as was pointed out by H.Gl\"ursing-L\"ur\ss en. We will see a counterexample in \S\ref{1dim}. She also showed that  one must consider a non-Noetherian extension of the ring $\mathbb{C}[\partial,\sigma]$ in order to obtain an Ehrenpreis-Malgrange type theorem for $D\Delta$-equations. Our interest, therefore, is to introduce a suitable ring extension $\mathcal{H}$ of $\mathbb{C}[\partial,\sigma]$ and prove Ehrenpreis-Malgrange type theorem for this ring. Our main result is the following:

%Let $\pi:\mathbb{R}^n\rightarrow \mathbb{R}^{n-1}$ be the projection which truncates the first coordinate. We can define the subsheaf $\mathcal{S}ol_{\mathbb{P},\mathscr{F}}$ of $\pi_*\mathscr{F}^{q\times 1}$ for $\mathscr{F}=C^\infty ,\mathcal{B}$ by the formula
%$$\mathcal{S}ol_{\mathbb{P},\mathscr{F}}= \Ker\left(\mathbb{P}:\pi_*\mathscr{F}^{q\times 1}\rightarrow\pi_*\mathscr{F}^{p\times 1}\right).$$

%\begin{cor}\label{solsheaf}
%For any open subset $\Omega^\prime\subset\mathbb{R}^{n-1}$, one has a canonical cohomology isomorphism
%$$H^{i}(\Omega^\prime ,\mathcal{S}ol_{\mathbb{P} ,\mathscr{F}})\simeq \Ext^i_{\mathcal{H}}(M, \pi_*\mathscr{F}(\Omega^\prime )),$$
%where $\mathscr{F}=C^\infty, \mathcal{B}$ and $i$ is any integer.
%\end{cor}

\begin{thm}\label{mainresult}
Let $\mathscr{F}$ be either $C^\infty$, $\mathcal{O}$, or $\mathcal{B}$ and  $\Omega$ be a convex subset of $\mathbb{R}^n$ (when $\mathscr{F}=C^\infty$, $\mathcal{B}$) or of $\mathbb{C}^n$ (when $\mathscr{F}=\mathcal{O}$) such that $\Omega +\mathbb{R}\mathbf{e}_1=\Omega$. Here, $\mathcal{B}$ is the sheaf of hyperfunctions. Let $\mathcal{H}$ denote a ring of generalized $D\Delta$-operators which will be introduced in $\S$\ref{ring}. Finally, let $\mathbf{E}$ denote an $\mathcal{H}$-module consisting of exponential polynomials. For any  matrix $\mathbb{P}\in M(r_1,r_0;\mathcal{H})$, we put $M=\mathcal{H}^{1\times p}/\mathcal{H}^{1\times q}\mathbb{P}.$ One has the following properties:
\begin{enumerate}
\item For any positive integer $i$, one has
\begin{equation}
\Ext^i_\mathcal{H}(M,\mathscr{F}(\Omega))=0.
\end{equation}
\item If $\mathscr{F}$ is either $C^\infty$ or $\mathcal{O}$, then $\Ker\left(\mathbb{P}(\partial ,\sigma ):\mathbf{E}^{r_0\times 1}\rightarrow\mathbf{E}^{r_1\times 1}\right)$ is dense in $\Ker(\mathbb{P}(\partial ,\sigma ):\mathscr{F}(\Omega)^{r_0\times 1}\rightarrow\mathscr{F}(\Omega)^{r_1\times 1}).$
\end{enumerate}
\end{thm}

Although Theorem \ref{EM} and Theorem \ref{mainresult} concern ring theoretic study of convolution equations, the technology of topological vector spaces enables us to reduce the proofs to solving a problem in complex harmonic analysis. Namely,  Ehrenpreis-Malgrange type theorem follows if, for any $r_1\times r_0$ matrix $\mathbb{P}$ with entries in $\mathbb{C}[\partial]$  (resp.~$\mathcal{H}$), one can show the identity
\begin{equation}
\mathcal{O}(\mathbb{C}^n)^{1\times r_1}\widehat{\mathbb{P}}\cap\mathcal{O}(\mathbb{C}^n)^{1\times r_0}_p=\mathcal{O}(\mathbb{C}^n)^{1\times r_1}_p\widehat{\mathbb{P}},
\end{equation}
$\text{where } \hspace{3mm}\widehat{ } \hspace{3mm}$ stands for the Fourier transform (of convolution operators)  and\\ $\mathcal{O}(\mathbb{C}^n)_p\subset\mathcal{O}(\mathbb{C}^n)$ is a set of holomorphic functions with a suitable growth encoded from the  Fourier transform.
This identity which is called ``division with bounds'', was developed by many mathematicians around 60's and gave rise to fruitful theorems in the theory of differential equations with constant coefficients. However, the classical proof of Theorem \ref{EM} requires a refined version of Noether's normalization lemma as well as \v{C}ech cohomology with bounds, neither of which can be adapted for the proof of Theorem \ref{mainresult}. The crucial step of the proof of Theorem \ref{mainresult} is to combine our ring theoretic study of $\mathcal{H}$ with integral representation technique developed by M. Andersson: the residue current method (\cite{AW}).

For the readers' convenience, we review a brief history of Ehrenpreis-Malgrange type theorem. Soon after the celebrated work of L. Ehrenpreis (\cite{Eh}), there naturally appeared some people investigating a generalization of Theorem \ref{EM} for more general convolution equations. Their interests were mainly about generalizing (2) of Theorem \ref{EM}. However, D. I. Gurevich gave the following counterexample against the density of exponential polynomial solutions for general convolution equations. 

\begin{prop}[\cite{Gr}]
There are two convolution operators $\mu_1$ and $\mu_2$ in  $\mathcal{E}^\prime$ so that the exponential polynomial solutions are not dense in the solution space $\{ f\in C^{\infty}(\mathbb{R}^n )\mid \mu_1*f=\mu_2*f=0\}$.
\end{prop}
\noindent
Due to this example, generalizations of  (2) of Theorem \ref{EM} have been limited to convolution equations whose Fourier transforms define complete intersection varieties and for a certain period of time, this constraint remained to be necessary. (In the case of discrete varieties, this constraint is significantly relaxed to the ``locally slowly decreasing condition'' in \cite{BKS}. See \cite[Definition 2.1]{BKS}.) The main technique of their  proofs was based on explicit integral formulae in the spirit of multidimensional residue theory: residue currents (cf. \cite{BGVY}, or \cite{BT}). In this direction, the most general version of such currents was finally introduced in \cite{AW}. The authors of \cite{AW} gave a new proof of (2) of Theorem \ref{EM}, which is a prototype of that of Theorem \ref{mainresult}.

On the other hand, if we restrict our attention to $D\Delta$-equations, ring theoretic study comes into play. It enables us to prove a variant of (1) of Theorem \ref{mainresult} for smaller function module such as $\mathscr{F}(\mathbb{C}^n)=\mathcal{O}(\mathbb{C}^n; \exp)$ (entire functions of exponential growth). See \cite[Theorem 7.4]{OB}. Though ring theoretic approach cannnot solve division with bounds in general, the observation of H. Gl\"using-L\"ur\ss en in \cite{GL} is of significant importance: in the theory of (ordinary) $D\Delta$-equations, one must consider a system of equations as a module over a non-Noetherian ring $\mathcal{H}_1$, a univariable version of $\mathcal{H}$. She also conducted an extensive study of ring theoretic properties of $\mathcal{H}_1$, which serves as a bridge between analysis of systems and homological algebra in this paper. This paper makes full use of these progresses in complex analysis and non-Noetherian ring theory and therefore, is positioned at an intersection of two different points of views.

Let us summarize the content of this paper. In Section 2 we introduce the ring $\mathcal{H}$ defined by Gl\"using-L\"ur\ss en and prove the integral representation of hyperfunction solutions as a superposition of exponential polynomial solutions in dimension 1. Section 2 is independent of other sections, but it would help the reader to understand the latter part of this paper. In Section 3 we study a generalization of $\mathcal{H}$ which only contains difference operators to one direction, and examine some of its cohomological properties. In Section 4 we recall the notion of the residue currents and give some proofs of necessary results. In Section 5 we combine the techniques of Sections 3 and 4 to obtain our version of the division with bounds.  On the way, we need an estimate of residue currents for exponential polynomials which was established in \cite{BY1}. At the beginning of Section 5, however, we will observe that division with bounds is impossible if there are more than two independent difference directions suggesting that our result is optimal. In Section 6 we prove the Theorem \ref{mainresult} and, as a corollary, a description of cohomology groups with values in solution sheaves of systems

The integral representation of Section 2 is similar to the result of Y. Okada (cf. \cite{Ok}), but the way of proving it is different. I also owe the idea of the proof of the integral representation of hyperfunction solutions to \cite{K} and \cite{O}.

\section{$D\Delta$-equations in dimension 1}\label{1dim}

In this section, we introduce a certain ring extension of the ring of $D\Delta$-operators when the number of independent variable is 1. We will also observe that various function spaces enjoy better cohomological properties over this ring than the usual ring of $D\Delta$-operators $\mathbb{C}[\partial,\sigma]$.
In concordance with the notation of H. Gl\"using-L\"ur\ss en, let us denote the standard coordinate of $\mathbb{C}$ by $z$ and consider the polynomial ring $ \mathbb{C}[z,\sigma ]$ of 2 variables.
We define an action of $ \mathbb{C}[z,\sigma ]$ on $C^\infty (\mathbb{R})$ by 
\begin{equation}
 (z\cdot f)(x)=\frac{df}{dx}(x),
\end{equation}
and
\begin{equation}
(\sigma\cdot f)(x)=f(x+1).
\end{equation}
As we remarked in the introduction, $C^\infty (\mathbb{R})$ is not injective over $ \mathbb{C}[z,\sigma ]$. The following simple counterexample is due to H. Gl\"using-L\"ur\ss en.

\begin{ex}
Take two matrices $\mathbb{P}=^{t}(\sigma -1,z),\mathbb{Q}=(z,1-\sigma )$ and consider an exact sequence
$$\mathbb{C}[z,\sigma]\stackrel{\times\mathbb{Q}}{\rightarrow }\mathbb{C}[z,\sigma]^{1\times 2}\stackrel{\times\mathbb{P}}{\rightarrow }\mathbb{C}[z,\sigma]\; (exact).$$
Applying $\Hom_{\mathbb{C}[z,\sigma]}(-,C^{\infty}(\mathbb{R}))$, it yields a complex
$$C^{\infty}(\mathbb{R})\stackrel{\mathbb{P}\cdot }{\rightarrow }C^{\infty}(\mathbb{R})^{2\times 1}\stackrel{\mathbb{Q}\cdot }{\rightarrow }C^{\infty}(\mathbb{R}) \;(not\; exact).$$
In fact, we have $^{t}(0,1)\in \Ker (\mathbb{Q}\cdot )\setminus \Image (\mathbb{P}\cdot )$.
\qed
\end{ex}

This example suggests that we need to introduce another appropriate ring of operators to recover the classical vanishing theorem of higher extension groups.  The above counterexample indicates that it should necessarily be a non-flat extension of the polynomial ring $\mathbb{C}[z,\sigma]$. Following H. Gl\"using-L\"ur\ss en, we introduce the ring of $D\Delta$-operators $\mathcal{H}$. For any element $q$ of the fraction field $\mathbb{C}(z,\sigma)$, we define the meromorphic function $q^*$ by $q^*(z)=q(z,e^z)$.
\begin{defn}
We define a ring extensions of $\mathbb{C}[z,\sigma]$ by
\begin{equation}
\mathcal{H}=\{ q=p\phi^{-1}\mid p\in\mathbb{C}[z,\sigma,\sigma^{-1}],\phi\in\mathbb{C}[z],q^*\in\mathcal{O}(\mathbb{C}_z)\}.
\end{equation}
\end{defn}

\begin{defn}
A commutative domain $R$ is called a B\'ezout domain if any finitely generated ideal of $R$ is principal.
\end{defn}

\begin{nonthm}[{\cite[Theorem 3.2.1]{GL}}]\label{GL}\hfill
\begin{enumerate}
\item $\mathcal{H}$ is a non-Noetherian ring but a B\'ezout domain.
\item $\mathcal{H}$ is an elementary divisor domain, i.e., for any $\mathbb{P}\in\mathcal{H}^{p\times q}$ one can find two matrices $V$ and $W$ such that
\begin{equation}
V\mathbb{P}W={\rm diag}(d_1,...,d_r,0,...,0),
\end{equation}
where $V$ and $W$ are products of elementary matrices and $d_i\in\mathcal{H}\setminus\{ 0\}$, $d_i|d_{i+1}$.
\end{enumerate}
\end{nonthm}

Roughly speaking, this theorem means that any system of $D\Delta$-equations can be transformed into a direct sum of single equations in a certain sense as we explain below.
Notice that function spaces $C^\infty (\mathbb{R}),$ $\mathcal{O}(D),$ $\mathcal{D}^\prime (\mathbb{R}),$ and $\mathcal{B}(\mathbb{R})$ are all $\mathcal{H}$-modules. Here, $D$ is a horizontal strip $D=\{z\in \mathbb{C}|a<\Image z<b\}$ $(a<b,a,b\in[-\infty ,\infty ])$.
The action of the ring $\mathcal{H}$ on these function spaces is defined as follows: when $q=p(z,\sigma )\phi(z)^{-1}\in\mathcal{H}$ and $f(x)$ is a function of one of four  kinds above, say $\mathcal{C}^\infty (\mathbb{R})$, we take $g\in\mathcal{C}^\infty (\mathbb{R})$ so that $\phi(z)\cdot g(x)=f(x)$. We define $q\cdot f=p(z,\sigma )\cdot g(x)$. Since $q^*$ is entire, this action is independent of the choice of $g$.
Hence, dealing with the ring $\mathcal{H}$ amounts to discussing generalized $D\Delta$- operators.

\begin{nonthm}\label{Theorem1Variable}
\hfill
\begin{enumerate}
\item (real case) 
Let F be either $C^\infty (\mathbb{R}) ,$ $\mathcal{D}^\prime (\mathbb{R}),$ or $\mathcal{B}(\mathbb{R})$. For any finitely presented $\mathcal{H}$-module $M$ and for any positive integer $i,$ one has the following  vanishing result
\begin{equation}
\Ext_{\mathcal{H}}^i(M,F)=0.
\end{equation}
\item (holomorphic case)
Let D be a subset of $\mathbb{C}$ defined by $D=\{z\in \mathbb{C}|a<\Image z<b\}$, where $a<b$, $a, b\in[-\infty ,\infty ])$. For any finitely presented $\mathcal{H}$-module $M$ and for any positive integer $i,$ one has the following  vanishing result
\begin{equation}
\Ext_{\mathcal{H}}^i(M,\mathcal{O}(D))=0.
\end{equation}
\end{enumerate}
\end{nonthm}

\vspace{3mm}

\begin{cor}[integral representation formula]\label{1stCorollary}
Take any $q\in\mathcal{H}$, and $f\in\mathcal{B}(\mathbb{R})$ and assume $q\cdot f=0$. If $\{(\alpha_k,m_k)\}_{k\geq 1}
\;(\alpha_k\in\mathbb{C},q^*(\alpha_k)=0,m_k=ord_{\alpha_k}q^* )$ are  zeros of $q^*$ with its multiplicities, then f has a representation
\begin{equation}
f=\sum_{k\geq 1}P_k(x)e^{\alpha_k\cdot x}\;\;(P_k(X)\in\mathbb{C}[X],degP_k<m_k).
\end{equation}
Here the sum is convergent in the sense of hyperfunctions.
\end{cor}
\vspace{-3mm}
%\begin{rem}
%If a subset of $\mathbb{C}$ is convex, connected, open, and $\mathbb{Z}$ invariant, it is $\{z\in \mathbb{C}|a<\Image z<b\}$ $(a<b,a,b\in[-\infty ,\infty ])$.
%\end{rem}

\vspace{1em}
\noindent
{\bf Proof of Theorem \ref{Theorem1Variable}}
We first prove (2).
By the definition of a finitely presented module, $M$ is of the form $M=\mathcal{H}^{1\times p}/\mathcal{H}^{1\times q}\cdot\mathbb{P}$ for some $\mathbb{P}\in\mathcal{H}^{q\times p}$. By Theorem \ref{GL}, we can assume $\mathbb{P}={\rm diag}(d_1,...,d_r,0,...,0)$.
Therefore, it is reduced to showing that for any $q\in\mathcal{H}\setminus \{0\},\;q:\mathcal{O}(D)\rightarrow\mathcal{O}(D)$ is surjective.
We can assume $q=p(z,\sigma )\cdot\sigma^l(l\in\mathbb{Z}),p=\sum_{j=0}^{M}p_j(z)\sigma^j$. 
The surjectivity of $q$ now follows from Theorem 1 of \cite{MS} (see also the remark after Theorem 1).

\noindent
The proof of (1) is parallel to that of (2). One can find the corresponding surjectivity results in \cite[Theorem 6.1.23.]{BG} (see also Proposition 3.1.31. of the same reference) for $F=C^\infty(\mathbb{R})$ or in \cite[Theorem 1]{Eh2} for $F=\mathcal{D}^\prime(\mathbb{R}).$ The surjectivity result for $F=\mathcal{B}(\mathbb{R})$ follows from the commutative diagram in the proof of Corollary \ref{1stCorollary}. 
\qed

\vspace{1em}

\noindent
{\bf Proof of Corollary \ref{1stCorollary}}
First, observe that $q:\mathcal{O}(\mathbb{R}\times \sqrt{-1}(\mathbb{R}\setminus 0))\rightarrow\mathcal{O}(\mathbb{R}\times \sqrt{-1}(\mathbb{R}\setminus 0))$ is surjective from Theorem \ref{Theorem1Variable} (2). Now let us consider the commutative diagram
\[
\xymatrix{
 & & & & \\
0\ar[r] & \Ker q \ar@{.>}[r] \ar[d] & \Ker q \ar@{.>}[r] \ar[d] & \Ker q  \ar[d]  & \\
0 \ar[r] & \mathcal{O}(\mathbb{C}) \ar[r] \ar[d]^q & \mathcal{O}(\mathbb{R}\times \sqrt{-1}(\mathbb{R}\setminus 0)) \ar[r] \ar[d]^q & \mathcal{B}(\mathbb{R}) \ar[r] \ar[d]^q & 0\\
0 \ar[r] & \mathcal{O}(\mathbb{C}) \ar[r] \ar[d] & \mathcal{O}(\mathbb{R}\times \sqrt{-1}(\mathbb{R}\setminus 0)) \ar[r] \ar[d] & \mathcal{B}(\mathbb{R}) \ar[r] \ar[d] & 0\\
& 0 \ar@{.>}[r] & 0 \ar@{.>}[r] & 0 & \\
 & & & & .\\
}
\]
The first two columns are exact and the two middle rows are also exact. By the snake lemma, we obtain an identity
\begin{equation}
\Ker(q:\mathcal{B}(\mathbb{R})\rightarrow \mathcal{B}(\mathbb{R}))\hspace{-.3em}=\hspace{-.3em}\frac{\Ker(q:\mathcal{O}(\mathbb{R}\times \sqrt{-1}(\mathbb{R}\setminus 0))\rightarrow \mathcal{O}(\mathbb{R}\times \sqrt{-1}(\mathbb{R}\setminus 0)))}{\Ker(q:\mathcal{O}(\mathbb{C})\rightarrow\mathcal{O}(\mathbb{C})},
\end{equation}
and the third column is exact. Applying Theorem 2.3 below to 
\begin{equation}
\Ker(q:\mathcal{O}(\mathbb{R}\times \sqrt{-1}(\mathbb{R}\setminus 0))\rightarrow \mathcal{O}(\mathbb{R}\times \sqrt{-1}(\mathbb{R}\setminus 0))),
\end{equation} one gets the existence of such series representation.
\qed

\begin{thm}[{\cite[Proposition 6.4.17]{BG}}, see also Remark 6.2.11]
Take any $\mu\in\mathcal{O}^\prime(\mathbb{C})$, and let $K\subset\mathbb{C}$ be a convex carrier of $\mu$ and $\Omega\subset\mathbb{C}$ be a convex open set. We assume $\widehat{\mu }$ is regularly decreasing. Denoting the zeros of $\widehat{\mu }$ with their multiplicities by $\{(\alpha_k,m_k)\}_{k\geq 1}\;
(\alpha_k\in\mathbb{C},\widehat{\mu }(\alpha_k)=0,m_k=ord_{\alpha_k}q^* )$, there is an increasing sequence  $1=k_1<k_2<\cdots$ such that for any $f\in\mathcal{O}(\Omega +K)$ with $\mu *f=0$, one has a unique representation 
\begin{equation}
f=\sum_{n\geq 1}\sum_{k_n\leq k <k_{n+1}}P_k(x)e^{\alpha_k\cdot x}\;\;\;(P_k(X)\in\mathbb{C}[X],degP_k<m_k).
\end{equation}
Furthermore, when $\widehat{\mu }$ satisfies slowly decreasing condition, the grouping of terms is not necessary, i.e., f has a unique  representation of the form
\begin{equation}
f=\sum_{k\geq 1}P_k(x)e^{\alpha_k\cdot x}\;\;\;(P_k(X)\in\mathbb{C}[X],degP_k<m_k).
\end{equation}
\qed
\end{thm}

We do not define the terminology "slowly decreasing" in this paper since the notion is not necessary in the following sections. We only note that any element of $\mathcal{H}$ is "slowly decreasing" in the sense of \cite[Definition 2.2.13]{BG}.

\section{Ring $\mathcal{H}_n$}\label{ring}

We introduce a ring of generalized partial $D\Delta$-operators in the spirit of H. Gl\"using-L\"ur\ss en and investigate the basic properties of this ring. The reader should be aware that most of the following propositions do not hold if we consider the smaller ring $\mathbb{C}[z_1,\cdots ,z_n,e^{z_1}]$.

\begin{defn}
We define the ring $\mathcal{H}_n$ of generalized partial $D\Delta$-operators inductively as follows:

\noindent
For $n=1$, we put $\mathcal{H}_1=\mathcal{H}$.

\noindent
For $n\geq 2$, we put $\mathcal{H}_n=\mathcal{H}_1[z_2,\cdots ,z_n]$.
\end{defn}
\begin{rem}
\indent
\begin{enumerate}
\item From now on, we use the symbol $\mathcal{H}$ for $\mathcal{H}_n$ when there is  no risk of confusion.
\item We embed $\mathcal{H}$ into $\mathcal{O}(\mathbb{C}^n)$ by replacing $\sigma$ by $e^{z_1}$.
\item As we shall see below, the coordinates $z$ correspond to partial derivatives and $\sigma $ corresponds to a difference operator to the direction of the first coordinate. Therefore, the ring $\mathcal{H}_n$ is a ring of generalized $D\Delta$-operators whose frequencies generate a rank one additive subgroup of $\mathbb{C}^n$.
\end{enumerate}
\end{rem}

\noindent
In order to state basic properties of the ring $\mathcal{H}$, we need some lemmas from algebra. The most fundamental property of $\mathcal{H}$ is the so called coherency.

\begin{defn}
Let $R$ be a commutative ring.
\begin{enumerate}
\item An $R$ module $M$ is said to be coherent if there is an exact sequence of the form $R^k\rightarrow R^l\rightarrow M\rightarrow 0$ where $k$ and $l$ are positive integers.
\item The ring $R$ is said to be coherent if for any positive integers $k$ and $l$, and for any morphism $f:R^k\rightarrow R^l$, $\Ker f$ is finitely generated.
\end{enumerate}
\end{defn}

\noindent
The following proposition is taken from \cite{G} and it ensures the coherency of $\mathcal{H}$.

\begin{prop}[{\cite[Corollary 7.3.4]{G}}]\label{coherent}
Let $R$ be a semihereditary ring and let $x_1,\ldots ,x_n$ be indeterminates. Then $R[x_1,\ldots ,x_n]$ is a coherent ring.
\end{prop}

\noindent
Recall that a commutative ring $R$ is semihereditary if every finitely generated ideal of $R$ is a projective $R$-module. For example, a B\'ezout domain is a typical example of a semihereditary ring.

\begin{nonprop}\label{firstprop}
\hfill
\begin{enumerate}
\item $\mathcal{H}$ is a coherent ring.
\item $\mathcal{H}_1\subset \mathcal{O}(\mathbb{C})$ is a flat ring extension.
\item $\mathcal{H}\subset\mathcal{O}(\mathbb{C}^n)$ is a flat ring extension.
\end{enumerate}
\end{nonprop}

\begin{prf}
(1) is immediate from Proposition \ref{coherent} and Theorem \ref{GL} (1).

\noindent
(2) goes as follows: First, we prove that for any finitely generated ideal $I$ of $\mathcal{H}_1$, $\Tor_1^{\mathcal{H}_1}(\mathcal{H}_1/I,\mathcal{O}(\mathbb{C}))=0$. Since $\mathcal{H}_1$ is B\'ezout, we can assume $I=(q)$ for some $q\in\mathcal{H}_1$. Consider the exact sequence
\begin{equation}
0\rightarrow \mathcal{H}_1\overset{q\times }{\rightarrow}\mathcal{H}_1\rightarrow\mathcal{H}_1/(q)\rightarrow 0.
\end{equation}
By tensoring $\mathcal{O}(\mathbb{C})$, we have a complex
\begin{equation}
0\rightarrow \mathcal{O}(\mathbb{C})\overset{q\times }{\rightarrow}\mathcal{O}(\mathbb{C})\rightarrow\mathcal{O}(\mathbb{C})/(q)\rightarrow 0.
\end{equation}
This is exact since $\mathcal{O}(\mathbb{C})$ is a domain. Now by a standard argument of homological algebra, one has the flatness since for any ideal $I$ of $\mathcal{H}_1$, $\mathcal{H}_1/I$ is an inductive limit of $\mathcal{H}_1$-modules of the form $\mathcal{H}_1/J$ where $J$ is a finitely generated ideal of $\mathcal{H}_1$.
%I think I can omit the following arguments. For a nonfinitely generated ideal $I$ of $\mathcal{H}_1$, $I$ can be expressed as a union $I=\displaystyle \bigcup_\alpha I_\alpha$ of finitely generated ideals $\{I_\alpha\}_\alpha$ where indices $\alpha$ run over a cofinal directed set with respect to inclusion. Therefore, we have an inductive system of coherent modules $\{\mathcal{H}_1/I_\alpha\}_\alpha$. We finally have \\
%$Tor_1^{\mathcal{H}_1 }(\mathcal{H}_1/I,\mathcal{O}(\mathbb{C}))=Tor_1^{\mathcal{H}_1 }(\underset{\underset{\alpha}{\rightarrow}}{lim}\mathcal{H}_1/I_\alpha ,\mathcal{O}(\mathbb{C}))=\underset{\underset{\alpha}{\rightarrow}}{lim}Tor_1^{\mathcal{H}_1 }(\mathcal{H}_1/I_\alpha,\mathcal{O}(\mathbb{C}))=0$.\\

\noindent
(3) We can observe that $\mathcal{H}=\mathcal{H}_1[z_2,\ldots ,z_n]\subset \mathcal{O}(\mathbb{C})[z_2,\ldots ,z_n]$ is a flat ring extension by (2). Therefore, it is enough to prove that $\mathcal{O}(\mathbb{C})[z_2,\cdots ,z_n]\subset\mathcal{O}(\mathbb{C}^n)$ is flat, namely
\begin{equation}
\Tor_1^{\mathcal{O}(\mathbb{C})[z_2,\ldots ,z_n]}(\mathcal{O}(\mathbb{C})[z_2,\ldots ,z_n]/I, \mathcal{O}(\mathbb{C}^n))=0
\end{equation}
 for all finitely generated ideals $I$ of $\mathcal{O}(\mathbb{C})[z_2,\ldots ,z_n]$. Since $\mathcal{O}(\mathbb{C})$ is a B\'ezout domain, $\mathcal{O}(\mathbb{C})[z_2,\ldots ,z_n]$ is a coherent ring by Proposition \ref{coherent}. This implies that for any  finitely generated ideal $I$ of  $\mathcal{O}(\mathbb{C})[z_2,\ldots ,z_n]$, $\mathcal{O}(\mathbb{C})[z_2,\ldots ,z_n]/I$ has a finite free resolution.
Now it is enough to prove that for a given exact sequence $L\rightarrow M\rightarrow N$ of coherent  $\mathcal{O}(\mathbb{C})[z_2,\ldots ,z_n]$ modules, 
\begin{equation}
\mathcal{O}(\mathbb{C}^n)\underset{\mathcal{O}(\mathbb{C})[z_2,\ldots ,z_n]}{\otimes}L\rightarrow   \mathcal{O}(\mathbb{C}^n)\underset{\mathcal{O}(\mathbb{C})[z_2,\ldots ,z_n]}{\otimes}M\rightarrow   \mathcal{O}(\mathbb{C}^n)\underset{\mathcal{O}(\mathbb{C})[z_2,\ldots ,z_n]}{\otimes}N
\end{equation}
is again exact.
The last complex is exact if and only if
\begin{equation}
{}_n\mathcal{O}_{z^0}\underset{\mathcal{O}(\mathbb{C})[z_2,\ldots ,z_n]}{\otimes}L\rightarrow  {}_n\mathcal{O}_{z^0}\underset{\mathcal{O}(\mathbb{C})[z_2,\ldots ,z_n]}{\otimes}M\rightarrow {}_n\mathcal{O}_{z^0}\underset{\mathcal{O}(\mathbb{C})[z_2,\ldots ,z_n]}{\otimes}N
\end{equation}
 for any $z^0\in\mathbb{C}^n$ since $\mathbb{C}^n$ is a Stein open set. See \cite[Result 6.1.7]{OB}.
Here, we can replace  ${}_n\mathcal{O}_{z^0}\underset{\mathcal{O}(\mathbb{C})[z_2,\ldots ,z_n]}{\otimes}-$ by $ {}_n\mathcal{O}_{z^0}\underset{ {}_1\mathcal{O}_{z_1^0}[z_2,\ldots ,z_n]}{\otimes}{}_1\mathcal{O}_{z_1^0}[z_2,\ldots ,z_n]\underset{ \mathcal{O}(\mathbb{C})[z_2,\ldots ,z_n]}{\otimes}-$, where  $z^0=(z_1^0,\ldots ,z_n^0)\in\mathbb{C}^n$. Without loss of generality, we can assume $z^0=0$.
Since $\mathcal{O}(\mathbb{C})\subset{}_1\mathcal{O}_0$ is a flat ring extension by \cite[Result 6.1.7]{OB}, we only need to consider the ring extension $ {}_1\mathcal{O}_{0}[z_2,\ldots ,z_n]\subset {}_n\mathcal{O}_{0}$ and prove its flatness.

We consider a triple $ {}_1\mathcal{O}_{0}[z_2,\ldots ,z_n]\subset {}_n\mathcal{O}_{0}\subset\mathbb{C}[[z_1,\ldots ,z_n]]$. If we denote the maximal ideal of $\mathbb{C}[[z_1,\ldots ,z_n]]$ by $\mathfrak{m}$, then we can observe that $\mathfrak{m}\cap{}_n\mathcal{O}_0$ and $\mathfrak{m}\cap{}_1\mathcal{O}_0[z_2,\ldots ,z_n]$ are maximal ideals of ${}_n\mathcal{O}_0$ and ${}_1\mathcal{O}_0[z_2,\ldots ,z_n]$, respectively.
By taking completions of these rings with respect to these maximal ideals, we see that  $ {}_n\mathcal{O}_{0}\subset\mathbb{C}[[z_1,\ldots ,z_n]]$ is faithfully flat and  $ {}_1\mathcal{O}_{0}[z_2,\ldots ,z_n]\subset \mathbb{C}[[z_1,\ldots ,z_n]]$ is flat. Therefore, we can conclude that $ {}_1\mathcal{O}_{0}[z_2,\ldots ,z_n]\subset {}_n\mathcal{O}_{0}$ is flat. \qed
\end{prf}

Notice that one can actually show that the ring extension $\mathcal{H}_1\subset \mathcal{O}(\mathbb{C})$ is a faithfully flat extension. See \cite[Theorem 7.7]{OB}. We do not know whether the ring extension $\mathcal{H}\subset\mathcal{O}(\mathbb{C}^n)$ is a faithfully flat ring extension. If it is, we can get a refined estimate of $proj.dim M$ for coherent $\mathcal{H}$-modules $M$.
\\
\\
%There is another proof of (3) of the above proposition\ref{firstprop} based on completion of coherent rings. 

%For reader's convenience, we give necessary preliminaries about completions of coherent rings which doesn't seem to be well-known.

%\begin{prop}[\cite{H}]\ 
%Let $R$ be a commutative coherent ring, $I$ be its finitely generated ideal, and $M$ be a coherent $R$ module. Then we have a canonical isomorphism
%$$\widehat{M}\simeq\widehat{R}\underset{R}{\otimes}M$$
%Here \; $\widehat{ }$ stands for the $I$-adic completion. In particular, $\widehat{R}\underset{R}{\otimes}-$ is an exact functor from the category of coherent $R$ modules to that of coherent $\widehat{R}$ modules. 
%\end{prop}

%(another proof of (3) of the above proposition\ref{firstprop})\\
%It is enough to confirm that $\mathcal{H}\subset\mathbb{C}$

The second property of the ring $\mathcal{H}$ is the following identity which shows that this ring is actually identical with the one discussed in Section 7 of \cite{OB}.

\begin{prop}
One always has the following identity:
\begin{equation}
\mathcal{H}=\mathcal{O}(\mathbb{C}^n)\cap \mathbb{C}(e^{z_1},z_1,\cdots , z_n).
\end{equation}
\end{prop}

\begin{prf}
Take any $f=\frac{P}{Q}\in\mathcal{O}(\mathbb{C}^n)\setminus \{ 0\}$ where $P,Q\in\mathbb{C}[e^{z_1},z_1,\cdots , z_n]$. By the generalization of the theorem of Ritt (cf. \cite[MAIN THEOREM]{BD}), we can assume $f=\frac{p}{q}$, for some $p=\sum_{j}a_j(z)e^{\langle\theta_j,z\rangle},\; a_j(z)\in\mathbb{C}[z_1,\dots,z_n],\; \theta_j\in\mathbb{C}^n,\; q\in\mathbb{C}[z_1,\ldots ,z_n]$. Furthermore, a careful reading of the proof of \cite[MAIN THEOREM]{BD} tells us that $\theta_j$ is obtained as a $\mathbb{Z}$-linear combination of frequencies of $P$ and $Q$ and their complex conjugates, meaning that we may assume $p\in\mathbb{C}[z_1,\dots,z_n,e^{z_1}]$.

We can assume $n=2$ since the proof of the proposition for $n\geq 3$ is essentially the same.  We want to show $\frac{\partial q}{\partial z_{2}}=0$. Assume the opposite, that is, assume $\frac{\partial q}{\partial z_{2}}\neq 0$. In this case, we should have $\frac{\partial p}{\partial z_{2}}\neq 0$ since the fraction $\frac{p}{q}$ is entire.

Now we factorize $p$ as $p=p_1(z_1,z_2) p_2(z_1,z_2,e^{z_1}),\; p_1\in\mathbb{C}[z_1,z_2],\; p_2\in\mathbb{C}[z_1,z_2,X]$ so that $$p_2=p_{21}^{m_1}\cdots p_{2l}^{m_l},\; \frac{\partial p_{2j}}{\partial X}\neq 0,\; p_{2j}\; \text{are irreducible in }\mathbb{C}[z_1,z_2,X],$$
and $q \;\text{and}\; p_1 \;\text{are relatively prime.}$ We also factorize $q$ as $q=q_{1}^{n_1}\cdots q_{k}^{n_k}$,  $q_i\in\mathbb{C}[z_1,z_2],$ and $q_{i}$ are irreducible. We further assume $\frac{\partial q_1}{\partial z_2}\neq 0$. Let $\Delta_{1} (z_1)\in\mathbb{C}[z_1]\setminus \{ 0\}$ be the discriminant of $q_{1}$ with respect to $z_2$. 
%Since it is clear that  $e^{z_1}\notin\overline{\mathbb{C}(z_1)}^{alg}$ (the bar stands for the algebraic closure), we have $\Delta_{2j}(z_1,e^{z_1})\not\equiv 0$.
 We can take $z_1^0\in\mathbb{C}$ so that 
$$\Delta_{1}(z_1^0)\neq 0.$$ 
Take $z_2^0\in\mathbb{C}$ so that $\frac{\partial q_1}{\partial z_2}(z_1^0,z_2^0)\neq 0$ and  $q_1(z_1^0,z_2^0)=0$.
We can solve the algebraic equation $q_1(z_1,z_2)=0$ on a neighbourhood of $z^0=(z_1^0,z_2^0)$ with respect to $z_2$, i.e., $z_2=z_2(z_1)$ around $z^0=(z_1^0,z_2^0)$ on the zero set $\{ q_1=0 \}$ of $q_1$. Since $p_1$ and $q_1$ are relatively prime, we have $p_1(z_1,z_2(z_1))\not\equiv 0$. Thus we have $p_2(z_1,z_2(z_1),e^{z_1})=0$ on $\{ q_1=0\}$ around $z^0$.
In particular, there is j so that $p_{2j}(z_1,z_2(z_1),e^{z_1})=0$ on $\{ q_1=0\}$ around $z^0$.\\
Here we need the following elementary lemma.
\end{prf}

\begin{lem}\ 
The function $z_2(z_1)$ can be analytically continued as a single valued function to a simply connected domain $D$ of $\mathbb{C}$ which is obtained by subtracting finitely many rays from the whole plane. Furthermore, $z_2(z_1)$ has polynomial growth at infinity, i.e., there are positive real numbers $C , N$ such that $|z_2(z_1)|\leq C(1+|z_1|)^N$.
\end{lem}

Now the identity $p_{2j}(z_1,z_2(z_1),e^{z_1})=0$ is still valid on the extended domain $D$. If we expand $p_{2j}$ as  $p_{2j}=\displaystyle\sum_{j=0}^{M}\beta_j(z_1,z_2)X^j$, this identity is expressed as $\displaystyle\sum_{j=0}^{M}\beta_j(z_1,z_2(z_1))e^{jz_1}=0$. We put $\alpha_j(z_1)=\beta_j(z_1,z_2(z_1))$, and claim that there is at least one $j\geq 1$ such that $\alpha_j(z_1)\not\equiv 0$. Suppose, conversely, that $\alpha_j\equiv 0$ for all $j\geq 1$, then one necessarily has $\alpha_0\equiv 0$. This implies that for all $j$, $q_1\underset{\mathbb{C}[z_1,z_2]}{|}\beta_j(z_1,z_2)$. This contradicts the assumption that $p_{2j}$ is irreducible. Thus, the claim was confirmed.

Now let $1\leq J\leq M$ be the largest number such that $\alpha_J\not\equiv 0$. Since  $\displaystyle\sum_{j=0}^{M}\alpha_j(z_1)e^{jz_1}=0$, we have $\displaystyle\sum_{j=0}^{J}\frac{\alpha_j(z_1)}{\alpha_J(z_1)}e^{(j-J)z_1}=0$. When $\Real z_1$ tends to $+\infty$, we get 0=1. This is a contradiction.\qed
\vspace{1em}

For a commutative ring $R$, let us denote by ${\rm gl.dim}R$ the global dimension of $R$ (cf. \cite{Wei} Definitions 4.1.1). Amongst all homological properties of $\mathcal{H}$, the most important one is the following.

\begin{thm}\label{gldim}
One always has the following inequality:
\begin{equation}
{\rm gl.dim}\mathcal{H}_n\leq n+1.
\end{equation}
\end{thm}

In order to prove Theorem \ref{gldim} above, it is enough to estimate ${\rm gl.dim}\mathcal{H}_1$ since we have the following general estimate.

\begin{thm}[{\cite[Theorem 1.3.16]{G}}]\label{ineq}
If $R$ is a commutative ring, and $T$ is an indeterminate, then
\begin{equation}
{\rm gl.dim} R[T]\leq {\rm gl.dim} R +1.
\end{equation}
\end{thm}
 
Furthermore, it is enough to estimate projective dimensions of $\mathcal{H}_1/I$ for all ideals $I$ of $\mathcal{H}_1$ (see \cite{Wei} Theorem 4.1.2). Note that the structure of finitely generated ideals are completely understood since $\mathcal{H}_1$ is B\'ezout.  The key point of the proof is that the structure of nonfinitely generated ideals was also deeply investigated by H. Gl\"using-L\"ur\ss en  (cf.\cite{GL}). In general, for any commutative ring $R$ the notation $a\underset{R}{|}b$ means that an element $b$ of $R$ is divisible by another element $a$ of $R$.

\begin{prop}[{\cite[Theorem 3.4.10]{GL}}]\label{prop34}
For any ideal $0\neq I\subset\mathcal{H}_1$, one can find a polynomial $p\in\mathbb{C}[z,\sigma]\setminus\{ 0\}$, and a set  $M\subset D_p\overset{def}{=}\{ \phi\in\mathbb{C}[z]\mid\phi :\text{ monic and } \phi\underset{\mathcal{H}_1}{|}p \}$ such that
\begin{enumerate}
\item $1\in M$, 
\item For any $\phi\in M$, and for any $\psi\in M,$ one has ${\rm LCM}(\phi ,\psi)\in M,$ 
\item For any $\phi\in M,$ and for any $\psi\in \mathbb{C}[z]\setminus\{ 0\}$ such that $\psi \underset{\mathbb{C}[z]}{|}\phi,$ one has $\psi\in M,$ 
\end{enumerate}
and one has the identity $I=\langle\langle p\rangle\rangle_{(M)}\overset{def}{=}\{ h\frac{p}{\phi }\mid h\in\mathcal{H}_1, \phi\in M\}$.
\end{prop}

\noindent
{\bf Proof of Theorem \ref{gldim}.}
In view of Theorem \ref{ineq}, it is enough to prove the inequality ${\rm gl.dim}\mathcal{H}_1\leq 2$. Let $I\subset\mathcal{H}_1$ be a non-zero ideal. By Proposition \ref{prop34}, we can assume $ I=\langle\langle p\rangle\rangle_{(M)}$ for some $p\in\mathbb{C}[z,\sigma]\setminus\{ 0\}$ and $M$ with properties (1), (2), and (3).

We consider an exact sequence
\begin{equation}\label{ExactSeq}
\mathcal{H}_1^{1\times M}\rightarrow\mathcal{H}_1\rightarrow\mathcal{H}_1/I\rightarrow 0\;(exact).
\end{equation}
Here the first morphism is given by $\displaystyle\sum_{\phi\in M}h_{\phi }e_{\phi }\mapsto\displaystyle\sum_{\phi\in M}h_\phi \frac{p}{\phi}$, where $\{ e_\phi\}_{\phi\in M}$ stands for a free basis of $\mathcal{H}_1^{1\times M}$. We put $K=\Ker(\mathcal{H}_1^{1\times M}\rightarrow\mathcal{H}_1)$ and for a finite subset $F$ of $M$, we also put  $\hat{\phi}=\displaystyle\prod_{\underset{\psi\in F}{\psi\neq\phi}}\psi$. 
We can observe 
\begin{equation}
\displaystyle\sum_{\phi\in F}h_\phi \frac{p}{\phi}=0\iff(\displaystyle\sum_{\phi\in F}h_\phi\hat{\phi}) \frac{p}{\displaystyle\prod_{\psi\in F}\psi}=0\iff \displaystyle\sum_{\phi\in F}h_\phi\hat{\phi}=0.
\end{equation}
We define a submodule of $\mathbb{C}[z]^{1\times F}$ by $K^{F,\;pol}=\left\{\displaystyle\sum_{\phi\in F}h_{\phi }e_{\phi }|\displaystyle\sum_{\phi\in F}h_{\phi }\hat{\phi }=0,\;h_\phi\in\mathbb{C}[z] \right\}$. Now we claim that the following identity is valid:
\begin{equation}
K^{F,\;pol}\otimes_{\mathbb{C}[z]}\mathcal{H}_1=K^F\overset{def}{=}\left\{\displaystyle\sum_{\phi\in F}h_{\phi }e_{\phi }\in\mathcal{H}_1^{1\times F}|\displaystyle\sum_{\phi\in F}h_{\phi }\hat{\phi }=0,\;h_\phi\in\mathcal{H}_1 \right\}.
\end{equation}
In fact, putting $F=\{ \phi_1,\ldots ,\phi_m\}$,
 we have an exact sequence
\begin{equation}
0\rightarrow K^{F,\; pol}\rightarrow\mathbb{C}[z]^{1\times F}
\xrightarrow{\times 
\begin{pmatrix}
\hat{\phi}_1\\
\vdots\\
\hat{\phi}_m
\end{pmatrix}
}\mathbb{C}[z]^{1\times 1}\;(exact).
\end{equation}
Taking into account that $\mathbb{C}[z]\subset\mathcal{H}_1$ is a flat ring extension, and tensoring $-\underset{\mathbb{C}[z]}{\otimes}\mathcal{H}_1$ to this sequence, we obtain the claim.

Let us go on with the proof of the theorem. Now the totality of finite subsets of $M$ forms a cofinal partially ordered set with respect to inclusion. So we have inductive systems $\{ K^F\}$ and $\{ K^{F,\; pol}\}$ where morphisms are natural inclusions $K^{F_1}\hookrightarrow K^{F_2}$ (resp. $K^{F_1,\; pol}\hookrightarrow K^{F_2,\; pol}$) for a pair of finite subsets $F_1\subset F_2$ of $M$.
Here we have the sequence of identities 
\begin{equation}
K=\underset{F\subset M;\; finite}{\bigcup }K^F=\underset{\underset{F\subset M;\; finite}{\longrightarrow}}{\lim}(K^{F,\; pol}\otimes_{\mathbb{C}[z]}\mathcal{H}_1)=\left(\underset{\underset{F\subset M;\; finite}{\longrightarrow}}{\lim}K^{F,\; pol}\right)\otimes_{\mathbb{C}[z]}\mathcal{H}_1.
\end{equation}
If we put $K^{pol}=\underset{\underset{F\subset M;\; finite}{\longrightarrow}}{\lim}K^{F,\; pol}$, since ${\rm gl.dim} \mathbb{C}[z]=1$, we have a projective resolution:
\begin{equation}
0\rightarrow P\rightarrow\mathbb{C}[z]^{1\times \Lambda }\rightarrow K^{pol}\rightarrow 0\; (exact),
\end{equation}
where $P$ is a projective $\mathbb{C}[z]$-module.
By tensoring $\mathcal{H}_1$, we have a projective resolution for $K$:
\begin{equation}
0\rightarrow P\otimes_{\mathbb{C}[z]}\mathcal{H}_1\rightarrow\mathcal{H}_1^{1\times \Lambda }\rightarrow K\rightarrow 0\; (exact),
\end{equation}
Attaching this resulting sequence to (\ref{ExactSeq}), we have a projective resolution of $\mathcal{H}_1/I$ of length 2:
\begin{equation}
0\rightarrow P\otimes_{\mathbb{C}[z]}\mathcal{H}_1\rightarrow\mathcal{H}_1^{1\times \Lambda }\rightarrow\mathcal{H}_1\rightarrow\mathcal{H}_1/I\rightarrow 0\;(exact).
\end{equation}
\qed

Combining Theorem \ref{gldim} with Theorem \ref{LequinSimis} below we get the following result.

\begin{thm}\label{thm53} 
Any coherent $\mathcal{H}$-module has a finite free resolution of length $\leq n+1$, that is, if $M$ is a coherent $\mathcal{H}$-module, we have a free resolution
\begin{equation}
0\rightarrow \mathcal{H}^{1\times r_{N}}\xrightarrow{\times \mathbb{P}_N}\cdots \rightarrow \mathcal{H}^{1\times r_1}\xrightarrow{\times \mathbb{P}_1}\mathcal{H}^{1\times r_0}\rightarrow M\rightarrow 0\;(exact),
\end{equation}
where $N\leq n+1$ and $r_j$ are positive integers.
\end{thm}

\begin{thm}[{\cite[Theorem b]{LS}}]\label{LequinSimis}
If $R$ is a B\'ezout domain, any finitely generated projective $R[X_1,\ldots ,X_n]$-module is free.
\end{thm}

We are now constructing the Hefer forms which admit a certain estimate of Paley-Wiener type. We can even determine their explicit forms as well so that the required estimate is satisfied in a trivial manner.

\begin{prop}\label{prop55}
For any given complex 
\begin{equation}
0\rightarrow \mathcal{H}^{1\times r_{N}}\xrightarrow{\times \mathbb{P}_N}\cdots \rightarrow \mathcal{H}^{1\times r_1}\xrightarrow{\times \mathbb{P}_1}\mathcal{H}^{1\times r_0}\rightarrow M\rightarrow 0,
\end{equation}
where $\mathbb{P}_j\in M(r_j,r_{j-1};\mathcal{H})$, there exist matrix-valued forms

$H^k_l(\zeta ,z)=\displaystyle\sum_{|I|=k-l}A_I d\zeta^I$ with $A_I\in M(r_k,r_l;\tilde{\mathcal{H} })$ such that 
\begin{align}
&H^k_l=0\; (k<l), \\
&H^l_l={\rm Id}_{r_l}, \\
&\delta_{\zeta -z}H^k_l=\mathbb{P}_k(\zeta )H^{k-1}_l-(-1)^{k-l-1 }H^k_{l+1}\mathbb{P}_{l+1}(z)\; (l< k),\label{TheEquation}
\end{align}
where $\delta_{\zeta -z}$ is the interior multiplication by $\sum_{i=1}^n(\zeta_i-z_i)\frac{\partial}{\partial\zeta_i}$. Here, we put $\tilde{\mathcal{H} }=A[\zeta_2,\cdots ,\zeta_n,z_2,\cdots ,z_n]$, where $A$ is the subring of $\mathcal{O}(\mathbb{C}^2_{(\zeta_1,z_1)})$ which is generated by $\{\frac{q(\zeta_1)-q(z_1)}{\zeta_1-z_1}|q\in\mathcal{H}_1\}$, $\{q(\zeta_1)|q\in\mathcal{H}_1\}$, and  $\{q(z_1)|q\in\mathcal{H}_1\}$.
\end{prop}

\begin{prf}
Fix $l$. We construct $H^k_l$ inductively. When $k=l+1$, the right hand side of the equation is just $\mathbb{P}_{l+1}(\zeta)-\mathbb{P}_{l+1}(z)$. It is enough to construct elements $h_1,\cdots, h_n$ of the ring $\tilde{\mathcal{H}}$ so that $\displaystyle\sum_{i=1}^nh_i(\zeta_i-z_i)=q(\zeta_1)\zeta^\alpha -q(z_1)z^\alpha$, where $q\in\mathcal{H}$ is any given element and $\alpha$ is a multiindex $\alpha=(\alpha_2,\cdots,\alpha_n)$. We assume $n=2$ since the essential part of the proof remains unchanged. In this case, we can take $h_1,\;h_2$ defined by
\begin{equation}
h_1=\left(\frac{q(\zeta_1)-q(z_1)}{\zeta_1-z_1}\right)z_2^\alpha,
\;\;\;h_2=q(\zeta_1)\left(\displaystyle\sum_{k=1}^{\alpha}\zeta_2^{\alpha-k}z_2^{k-1}\right).
\end{equation}
Let $k>l+1$ and consider general $n$. Since for each $i$, 
\begin{equation}
\tilde{\mathcal{H}}/((\zeta_n-z_n),\cdots ,(\zeta_{i+1}-z_{i+1}))\simeq A[\zeta_2,\cdots ,\zeta_n, z_2,\cdots ,z_i]
\end{equation}
is a domain, $\{ (\zeta_n -z_n),\cdots ,(\zeta_1-z_1)\}$ is a regular sequence.
By a general result on Koszul complex we have an exact sequence:
\begin{equation}
0\rightarrow\tilde{\mathcal{H}}d\zeta_1\wedge\cdots\wedge d\zeta_n\overset{\delta_{\zeta -z}}{\rightarrow}\displaystyle\sum_{|I|=n-1}\tilde{\mathcal{H}}d\zeta^I\overset{\delta_{\zeta -z}}{\rightarrow}\cdots\overset{\delta_{\zeta -z}}{\rightarrow}\displaystyle\sum_{i=1}^n\tilde{\mathcal{H}}d\zeta_i\overset{\delta_{\zeta -z}}{\rightarrow}\tilde{\mathcal{H}}\; (exact).
\end{equation}
Since the right hand side of the equation (\ref{TheEquation}) is $\delta_{\zeta -z}$ closed, we can take the desired matrix valued form $H^k_l(\zeta ,z)=\displaystyle\sum_{|I|=k-l}A_I d\zeta^I$.\qed
\end{prf}

\section{Residue currents and integral formulae}

In the proof of our main theorem, we make full use of the theory of residue currents. To begin with, we need to review the concrete construction of residue currents following \cite{AW}. From now on, $X$ always denotes a connected complex manifold of dimension $n$. Let  $E$ and $Q$ be two holomorphic hermitian vector bundles on $X$. We denote by $(\cdot,\cdot)$ the hermitian metric on $E$. Take any morphism $f\in \Hom(E,Q)$. Let $\sigma :Q\rightarrow E$ be the minimal inverse of $f$, i. e., $\sigma\xi$ is the minimal solution $\eta$ of $f\eta=\xi$ if $\xi\in\Image f$ and $\sigma\xi=0$ if $\xi$ is orthogonal to $\Image f.$ The minimal inverse $\sigma$ is divergent  when the rank of $f$ degenerates. The order of divergence can be described in terms of the determinant of $f$. Let us remember that the optimal rank $\rho =\underset{x\in X}{\sup}\rank f_x$ of $f$ is well-defined. We can also define the canonical section $F=\frac{1}{\rho !}\wedge^\rho f=det^\rho f$ of $\wedge^\rho E^*\otimes \wedge^\rho Q$. We put $Z=\{ x\in X|\wedge^\rho f_x=0\}$. Since $Z$ is a locally common zero set of finitely many holomorphic functions, it is a proper analytic subset of $X$ and $\sigma $ is smooth outside $Z$. Recall that for a given section $s$ of a hermitian vector bundle $E$, one can define a section $s^*$ of $E^*$ by $s^*=(\cdot ,s)$; it is called the dual section of $s$. The following lemma is the key tool for obtaining a good estimate of residue currents as we will see in Section 5.

\begin{lem}[{\cite[Lemma 4.1]{Andersson}}]\label{TheLemma}
Let $s$ be a global section of $\Hom(Q,E)\simeq E\otimes Q^*$ such that $f\circ s|_{\Image f}=|F|^2Id_{\Image f}$ and $s|_{(\Image f)^\perp}=0$ with pointwise minimum norms on $X\setminus Z$, and let $S$ be a global section of $\wedge^\rho E\otimes\wedge^\rho Q^*$ such that $FS=|F|^2$ with point wise minimum norms on $X\setminus Z$. Then  one has the identity 
\begin{equation}
s=|F|^2\sigma
\end{equation}
on $X\setminus Z$ and both $s$ and $S$ are smooth across $Z$.
\end{lem}

Though we do not include the proof of Lemma \ref{TheLemma}, we recall the explicit formula for $s$ which can be found in \cite{Andersson}. Take any smooth local frame $\{ \epsilon_j\}_{j=1}^r\; (r=\rank Q)$ of $Q$. With the aid of this frame, we can write $f$ locally as 
\begin{equation}
f=\displaystyle\sum_{j=1}^\rho (f^j)\otimes \epsilon_j,
\end{equation}
where $f^j\in E^*$. If we denote by $\delta_{f^j}$ (resp. $\delta_{\epsilon_j}$) the interior multiplication by $f^j$ (resp. $\epsilon_j$), the formula of the section $s$ is given by
\begin{equation}
s=\left(\displaystyle\sum_{j=1}^r\delta_{f^j}\otimes\delta_{\epsilon_j}\right)^{\rho -1}S/\rho !.
\end{equation}
Notice that the operator $\displaystyle\sum_{j=1}^r\delta_{f^j}\otimes\delta_{\epsilon_j}$ does not depend on a particular choice of the local representation $f=\displaystyle\sum_{j=1}^r f^j\otimes\epsilon_j$ even if local frame $\{ \epsilon_j\}_j$ of $Q$ is non-holomorphic. We also remark that $S$ is nothing but the dual section $F^*$ of $F$.

Before discussing the residue currents, we introduce a convention on composition of vector bundle-valued forms. Let us denote by $\mathcal{D}^\prime_k$ the sheaf of $k$-currents on X. For any given three vector bundles $E_1,\; E_2,\; E_3$, for any currents $\omega\in\mathcal{D}^\prime_k$, $\eta\in\mathcal{D}^\prime_l$, and for any morphisms $f\in \Hom(E_1,E_2)$ and $g\in \Hom(E_2,E_3)$, we take the following composition rule:
$$(\omega\otimes f)\cdot(\eta\otimes g)=\eta\wedge\omega\otimes f\circ g.$$

With these preparations, we can give the construction of residue currents following M. Andersson and E. Wulcan. 
In \cite{AW} they constructed the residue current by using the terminologies of super-connection (cf. \cite{Q}). We do not need to employ this formalism in this paper. Instead, we only give the explicit way of constructing it.

Let us consider a generically exact complex of hermitian vector bundles
\begin{equation}\label{NonSurjSeq}
0\rightarrow E_N\stackrel{f_N}{\rightarrow}\cdots \rightarrow E_1\stackrel{f_1}{\rightarrow}E_0.
\end{equation}
%\stackrel{f_0}{\rightarrow}E_{-1}\stackrel{f_{-M+1}}{\rightarrow}E_{-M}\rightarrow 0
Take the minimal inverse $\sigma_k$ of each morphism $f_k:E_{k}\rightarrow E_{k-1}$, and let $Z$ be the set of points of $X$ where some $f_k$ do not have optimal rank.
We define a smooth $\Hom(E_l,E_k)$-valued $(0,k-l-1)$-form $u_k^l$ for $(k\geq l+1)$ on $X\setminus Z$ by the following formula:
\begin{equation}
u_k^l=(\bar{\partial }\sigma_k)\cdots (\bar{\partial }\sigma_{l+2})\sigma_{l+1}.
\end{equation}

\noindent
With this notation, we have

\begin{prop}[\cite{AW}]\label{prop42}
For any local holomorphic function $F\in\mathcal{O}_X$ such that $Z\subset\{ F=0\}$ holds locally, $|F|^{2\lambda }u_k^l$ and $\bar{\partial }|F|^{2\lambda }\wedge u_k^l$ can be continued analytically to a current across $Z$ with a holomorphic parameter in $\Real\lambda >-\varepsilon$ for some small $\varepsilon >0$. The value at $\lambda=0$ of these currents are independent of the particular choice of $F$.
Furthermore, if we put
\begin{equation}
U_k^l=|F|^{2\lambda }u_k^l|_{\lambda =0}\;\; \text{for}\; k\geq l+1,
\end{equation}
 and 
\begin{equation}
R^l_k=\bar{\partial }|F|^{2\lambda }\wedge u_k^l|_{\lambda =0}\;\; \text{for}\; k\geq l+1,
\end{equation}
if the holomorphic function $F$ above can be factorized as $F=\tilde{F}_1\cdots \tilde{F}_p$, if $t_1,\cdots ,t_p$ are positive real numbers, and if $|\tilde{F}|^{*(t\lambda)}$ denotes $|\tilde{F}_1|^{t_1\lambda }\cdots |\tilde{F}_p|^{t_p\lambda}$, then $|\tilde{F}|^{*(t\lambda)}u_k^l$ and $\bar{\partial }|\tilde{F}|^{*(t\lambda)}\wedge u_k^l$ can be continued analytically to a distribution across $Z$ with a holomorphic parameter in $\Real\lambda >-\varepsilon$ for some $\varepsilon >0$ and one has the formulae
\begin{equation}\label{Unusual1}
U_k^l=|\tilde{F}|^{*(t\lambda)}u_k^l|_{\lambda =0}\;\; \text{for}\; k\geq l+1
\end{equation}
and
\begin{equation}\label{Unusual2}
R^l_k=\bar{\partial }(|\tilde{F}|^{*(t\lambda)})\wedge u_k^l|_{\lambda =0}\;\; \text{for}\; k\geq l+1.
\end{equation}
\end{prop}

\vspace{1em}
\noindent
The formulae (\ref{Unusual1}) and (\ref{Unusual2}) are not explicitly written in \cite{AW}. However, a careful reading of the constructions of $U_k^l$ and $R^l_k$ in \cite{AW} tells us that such formulae are valid. We put 
\begin{equation}
U=\displaystyle\sum_{l\geq 0}\sum_{k\geq l+1}U_k^l
\end{equation}
and 
\begin{equation}
R=\displaystyle\sum_{l\geq 0}\sum_{k\geq l+1}R_k^l,
\end{equation}
and call them associated currents to the complex (\ref{NonSurjSeq}). In particular, the current $R$ is called the residue current of the complex. The following result states that $R$ measures the exactness (hence non-exactness) of the complex (\ref{NonSurjSeq}).

\begin{thm}[{\cite[Theorem 1.1]{AW}}]\label{thm41}
The complex (\ref{NonSurjSeq}) is exact on $X$ if and only if $R_k^l=0$ for all positive integers $k>l$.
Furthermore, if the complex $(\ref{NonSurjSeq})$ is exact, a holomorphic section $\phi\in E_0$ is in $\Image f_1$ if and only if  $\phi$ is generically in $\Image f_1$ and for all $k>0,$ one has $R^0_k\phi =0$.
\end{thm}

%Here, notice that $R^0_k\phi$ stands for a current, not the evaluation of the current $R^0_k$ at $\phi$. 
Thus, the residue current $R$ actually corresponds to the notion of Noetherian operators defined originally by L. Ehrenpreis (see \cite{P}, and \cite{Eh} ).
%In the statement of the theorem above, that $\phi\in E_0$ is in $\Image f_1$ means $\phi$ is represented locally as $\phi =f_1\psi $, where $\psi\in E_1$ is a local holomorphic section.
This theorem gives us an abstract criterion of membership, but we need to know explicitly what $\phi$ in Theorem \ref{thm41} is in view of applications to the theory of equations. To this end, we use the integral representation technique that was introduced by B. Berndtsson and M. Andersson and further developed by M. Andersson (cf. \cite{AI}, \cite{AII}).\\

Let $D\subset\mathbb{C}^n$ be a domain and we regard $z\in D$ as a parameter. We put 
\begin{align}
\mathcal{L}^m&=\displaystyle\bigoplus_{k\geq 0}\mathcal{D}^\prime_{(k,k+m)},\\
\nabla_{\zeta -z}&=\delta_{\zeta -z}-\bar{\partial}.
\end{align}
Here, $\delta_{\zeta-z}$ is the operator introduced in Proposition \ref{prop55}. Note that any element $\omega\in\mathcal{L}^m$ can be decomposed as $\omega=\omega_{0,m}+\omega_{1,m+1}+\cdots$, where $\omega_{k,k+m}\in\mathcal{D}^\prime_{(k,k+m)}$.

\begin{prop}[{\cite[Proposition 2.1]{AII}}]\label{intrep}
Let $z\in D$ be a fixed point and $g=g_{0,0}+\cdots +g_{n,n}\in\mathcal{L}^0(D)$ be a current with compact support in D such that $\nabla_{\zeta -z} g=0$, $g$ is smooth around z, and $g_{0,0}(z)=1$. In this setting, for any holomorphic function $\phi\in\mathcal{O}(D)$, one has an integral representation formula
\begin{equation}
\phi (z)=\int_{D}g\phi =\int_{D}g_{n,n}\phi .
\end{equation}
\end{prop}

\noindent
Any $g\in\mathcal{L}^{0}(D)$ with properties in the proposition above is called a weight with respect to $z$. 
%Bochner-Martineli type weight has a singularity at $\zeta =z$. It is more convenient to use a non-singular weight.  
In \cite{AI} and \cite{AII} it was explained that one can obtain various classical integral formulae in a unified manner thanks to the  proposition above.

Let us now introduce yet another tool to get an explicit solution of the membership problem. We consider a generically exact complex.
\begin{equation}\label{SurjSeq}
0\rightarrow E_N\stackrel{f_N}{\rightarrow}\cdots \rightarrow E_1\stackrel{f_1}{\rightarrow}E_0\rightarrow 0.
\end{equation}
The readers should be aware that this complex is different from the complex (\ref{NonSurjSeq}). 
We recall the following general existence theorem of Hefer forms.

\begin{prop}[{\cite[Proposition 5.3]{AII}}]\label{prop44}
Assume $D$ is Stein and consider either (\ref{NonSurjSeq}) or (\ref{SurjSeq}). Then, for any integers $k$ and $l$, one can find $H^k_l(\zeta ,z)\in\mathcal{E}_{(k-l,0)}(\Hom(E_k,E_l))(D_\zeta)$ such that
$$
\begin{aligned}
&H^k_l(\zeta ,z) \text{ is holomorphic both in}\; \zeta \;\text{and}\; z.\\
&H^k_l(\zeta ,z)=0\; (k<l) \;\text{and} \;H^l_l(\zeta ,z)=Id_{E_l}.\\
&\delta_{\zeta -z}H^k_l(\zeta ,z)=H^{k-1}_l(\zeta ,z)f_k(\zeta )-(-1)^{k-l-1}f_{l+1}(z)H^k_{l+1}\;\; (l<k).
\end{aligned}
$$
\end{prop}

We define some important currents as follows:
$$
\begin{aligned}
&HU=\displaystyle\sum_{l}H_{l+1}U=\displaystyle\sum_{k\geq l+1}H^k_{l+1}U^l_k ,&  &H^k_{l+1}U^l_k\in\mathcal{D}^\prime_{(k-l-1,k-l-1)}(Hom(E_l,E_{l+1}))\\
&HR=\displaystyle\sum_{l}H_{l}R=\displaystyle\sum_{k\geq l+1}H^k_{l}R^l_k ,& &H^k_{l}R^l_k\in\mathcal{D}^\prime_{(k-l,k-l)}(End(E_l))\\
&f=\displaystyle\sum_{j\geq 1}f_j ,& &g^\prime (\zeta ,z)=f(z)HU+HUf+HR.
\end{aligned}
$$
Now we can check by a direct computation that $\nabla_{\zeta -z}g^\prime =0$. 
Furthermore, for any $z\in D\setminus Z,$ one has the identity $g^\prime_{0,0}(z,z)=\displaystyle\sum_{j\geq 0}{\rm Id}_{E_j}$, where $Z$ is the singular locus of currents $U$ and $R$. With this notation Proposition \ref{intrep} and Proposition \ref{prop44} lead to

\begin{prop}[{\cite[Proposition 5.4]{AII}}]\label{prop45}
Let $D$ be a Stein open subset of $\mathbb{C}^n$ and consider the generically exact complex (\ref{SurjSeq}).
\begin{enumerate}
\item If $g\in\mathcal{L}^0(D_\zeta )$ is a smooth weight with respect to $z\notin Z$, then for any holomorphic section $\phi\in\mathcal{O}(D_\zeta ,E_l),$ one has an integral representation
\begin{equation}\label{DivisionFormula}
\phi (z)=f_{l+1}(z)\int_{D_\zeta }H_{l+1}U\phi\wedge g + \int_{D_\zeta }H_lUf_l\phi\wedge g + \int_{D_\zeta}H_lR\phi\wedge g,
\end{equation}
\item If $g(\zeta ,z)\in\mathcal{L}^0(D_\zeta )$ is a smooth weight with respect to $z\in D$, $g $ is holomorphic in $z$, and if for all $z\in D$, $supp \{ g(\cdot ,z)\}$ is compact in $D$, then the (\ref{DivisionFormula}) is valid across $Z$.
\end{enumerate}
\end{prop}

\begin{prf}
Part (1) follows immediately from Proposition \ref{intrep} since $g^\prime\wedge g$ is a weight. To see Part (2) notice that in this case, both sides of the equation (\ref{DivisionFormula}) are holomorphic so we have the equation even across $Z$ by the  identity theorem. \qed
\end{prf}

Remember that we assumed $f_1$ is generically surjective in the argument above.  However, this assumption is not necessary. We shall see that we anyway can get a representation of $\phi$ as soon as it belongs to the image of $f_1$. To this end, we again  consider the complex (\ref{NonSurjSeq}).
We can extend this to another generically exact complex in a neighbourhood of each point by the following proposition.\\
\begin{prop}\label{prop45}
Let $X$ be a Stein manifold, and $K$ be a holomorphically convex compact subset of $X$. If $\mathcal{M}$ is a coherent $\mathcal{O}_X$-module on a neighbourhood of $K$, it can be embedded into a generically exact complex of the folowing type in a neighbourhood of $K$:
\begin{equation}
0\rightarrow\mathcal{M}\rightarrow\mathcal{O}_X^{1\times r_0}\rightarrow\cdots\rightarrow\mathcal{O}_X^{1\times r_{M-1}}\rightarrow\mathcal{E}\rightarrow 0,
\end{equation}
where $\mathcal{E}$ is a locally free sheaf in a neighbourhood of $K$.
\end{prop}

\begin{prf}
We put $\mathcal{N}^*=\mathcal{H}om(\mathcal{N},\mathcal{O}_X)$ for any $\mathcal{O}_X$-module $\mathcal{N}$.
By Theorem 7.2.1 in \cite{Ho}, one can take the following resolution of $\mathcal{M}$ in a neighbourhood of $K$:
\begin{equation}
\mathcal{O}_X^{1\times r_M}\rightarrow\cdots\rightarrow\mathcal{O}_X^{1\times r_{0}}\rightarrow\mathcal{M}^*\rightarrow 0\;\; (exact).
\end{equation}
Here, we put $M=\dim X.$ By Hilbert's syzygy theorem, we can conclude  $\mathcal{K}=Ker(\mathcal{O}_X^{1\times r_M}\rightarrow\mathcal{O}_X^{1\times r_{M-1}})$ is locally free.
By putting $\mathcal{E}=\mathcal{K}^*$, taking $(-)^*$ of the above complex, and composing the canonical map $\mathcal{M}\rightarrow\mathcal{M}^{**}$, we have a complex
\begin{equation}
0\rightarrow\mathcal{M}\rightarrow\mathcal{O}_X^{1\times r_0}\rightarrow\cdots\rightarrow\mathcal{O}_X^{1\times r_{M-1}}\rightarrow\mathcal{E}\rightarrow 0.
\end{equation}
Since this complex is exact where $\mathcal{M}$ is locally free, this is the desired complex.\qed
\end{prf}

With the help of this proposition, we can prove\\

\begin{prop}\label{prop46}
Let $D$ be a Stein open subset of $\mathbb{C}^n$ and $0\rightarrow E_N\stackrel{f_N}{\rightarrow}\cdots \rightarrow E_1\stackrel{f_1}{\rightarrow}E_0$ be a generically exact complex on $D$.\\
(1) If $g\in\mathcal{L}^0(D_\zeta )$ is a smooth weight with respect to $z\notin Z$, then, for any holomophic section $\phi\in\mathcal{O}(D_\zeta ,\Image f_1),$ one has an integral representation
\begin{equation}
\phi (z)=f_{1}(z)\int_{D_\zeta }H_{1}U\phi\wedge g.
\end{equation}
(2) If $g(\zeta ,z)\in\mathcal{L}^0(D_\zeta )$ is a smooth weight with respect to $z\in D$, $g $ is holomorphic in $z$, and if for all $z\in D$, $supp \{ g(\cdot ,z)\}$ is compact in $D$, then the formula of (1) is valid across $Z$.
\end{prop}

\begin{prf}
We first, embed the complex (\ref{NonSurjSeq}) into a generically exact sequence
\begin{equation}
0\rightarrow E_N\xrightarrow{f_N}\cdots \rightarrow E_1\xrightarrow{f_1}E_0\xrightarrow{f_0}\cdots\rightarrow E_{-n}\rightarrow 0
\end{equation}
on a neighbourhood of $\supp g$ by Proposition \ref{prop45}.
We can now prolong the Hefer forms taken above to $\{ H^k_l\}_{-n\leq l\leq k\leq N}$, where $H^{k}_l\in\mathcal{E}_{(k-l,0)}(\Hom(E_k,E_l))$ satisfy the relations of Proposition \ref{prop44}. Applying Proposition \ref{prop45} to this complex and $\phi\in\mathcal{O}(D_\zeta ,\Image f_1)$, we obtain the proposition.\qed
\end{prf}

\section{Division with bounds}\label{DWB}

In this section, we solve a certain kind of the division with bounds which naturally arises in the theory of $D\Delta$-equations. Before starting a discussion, we would like to note that the division with bounds is no longer possible if there are more than two independent frequencies. 
\begin{ex}
For any positive integer $n,$ let us denote by $\mathcal{O}_p(\mathbb{C}^n)$ the space of entire functions on $\mathbb{C}^n$ for which there exist constants $C>0$, $M>0$, and $k\in\mathbb{Z}_{>0}$  such that for every $\zeta\in\mathbb{C}^n$ one has the estimate
\begin{equation}
|f(\zeta)|\leq C(1+|\zeta|)^k\exp(M|\Image \zeta|).
\end{equation}

\noindent
Consider now an ideal $I$ of $\mathcal{O}_p(\mathbb{C}^2)$ generated by three elements $\frac{\sin\zeta_1}{\zeta_1}$, $\sin\zeta_2$, and $\zeta_2-\alpha\zeta_1$, where $\alpha\in\mathbb{R}\setminus\bar{\mathbb{Q}}$ is a Liouville number. We can easily check that 1 is in $(\mathcal{O}(\mathbb{C}^2)\cdot I)\cap\mathcal{O}_p(\mathbb{C}^2)$ since the zero locus of $I$ is empty. On the other hand, 1 can never belong to the ideal $I$. Otherwise we have a representation $1=f_1(\zeta_1,\zeta_2)\frac{\sin\zeta_1}{\zeta_1}+f_2(\zeta_1,\zeta_2)\sin\zeta_2+f_3(\zeta_1,\zeta_2)(\zeta_2-\alpha\zeta_1)$ for some $f_1,f_2,f_3\in\mathcal{O}_p(\mathbb{C}^2)$.  Substituting $\zeta_2=\alpha\zeta_1$, we can conclude that 1 belongs to the ideal $I^\prime$ of $\mathcal{O}_p(\mathbb{C})$ generated by $\frac{\sin\zeta_1}{\zeta_1}$ and $\sin\alpha\zeta_1$, where the definition of $\mathcal{O}_p(\mathbb{C})$ is similar to that of $\mathcal{O}_p(\mathbb{C}^2)$. However, L. Ehrenpreis showed that $I^\prime$ does not contain $1$, meaning that $I \subsetneq(\mathcal{O}(\mathbb{C}^2)\cdot I)\cap\mathcal{O}_p(\mathbb{C}^2)$ (cf. \cite[pp 319-320]{Eh}).
\end{ex}

This example suggests that the division with bounds is impossible when there are more than two independent frequencies. In spite of this sort of example, the division with bounds is possible for $D\Delta$-equations with one independent frequency as we shall see in this section. We begin with some elementary estimates related to the Hefer forms.

\begin{lem}\label{prop61}
Take any $q(z )\in\mathcal{H}_1$ and put $p(\zeta ,z)=\frac{q(\zeta )-q(z)}{\zeta -z}$. Assume that $q$ does not contain any negative power of $e^z$. Then one has the following estimates:
\begin{enumerate}
\item Regard $q\in\mathbb{C}(z)[e^z ]$. If $N$ denotes the degree of $q$ as a polynomial of $e^z$ with coefficients in $\mathbb{C}(z)$, then there are a non-negative integer $M\in \mathbb{Z}_{\geq 0}$ and a positive constant $ C>0$ such that for all $z\in\mathbb{C},$
\begin{equation}\label{BasicEstimate}
|q(z)|\leq C(1+|z|)^Me^{N|\Real z|}.
\end{equation}
\item With the same notation as above, there are a non-negative integer  $M\in \mathbb{Z}_{\geq 0}$ and a constant $ C>0$ so that the following inequality is valid for all $\zeta , z\in\mathbb{C}$:
\begin{equation}
|p(\zeta ,z)|\leq C(1+|\zeta |)^M(1+|z|)^Me^{N|Re \zeta |}e^{N|\Real z|}.
\end{equation}
\item For any non-negative integer $k\in\mathbb{Z}_{\geq 0}$, one can find a non-negative integer $M\in \mathbb{Z}_{\geq 0}$ and a constant $ C>0$ so that the following inequality is true for all $\zeta , z\in\mathbb{C}$:
\begin{equation}
\left|\frac{\partial^k}{\partial\zeta^k}p(\zeta ,z)\right|\leq C(1+|\zeta |)^M(1+|z|)^Me^{N|Re \zeta |}e^{N|\Real z|}.
\end{equation}
\end{enumerate}
\end{lem}

\begin{prf}
We omit the proof of (1) since it is straightforward. To see (2) first recall that we have a formula
\begin{equation}
q(\zeta )-q(z)=\left(\int_0^1q^{\prime }(z+t(\zeta -z))dt\right)(\zeta -z)
\end{equation}
by integration by parts. Note that $q^\prime (z)$ satisfies the estimate of (\ref{BasicEstimate}). Therefore, we have
\begin{align}
 &|q^{\prime }(z+t(\zeta -z))|&\leq&C(1+|z+t(\zeta -z)|)^Me^{N|\Real z+t \Real(\zeta -z)|}\\
 & &\leq& C(1+(1-t)|z|+t|\zeta |)^Me^{N|(1-t)\Real z+t \Real \zeta|}\\
 & &\leq& C(1+|\zeta |)^M(1+|z|)^Me^{N|\Real \zeta |}e^{N|\Real z|}. 
\end{align}
As for (3), if we prove the estimate when $k=1$, we inductively have the estimate for all $k$. When $k=1$, the inequality is immediately verified by means of Cauchy's integral formula as follows: 
\begin{align}
\left|\frac{\partial }{\partial\zeta }p(\zeta ,z)\right|&\leq C(1+|z|)^Me^{N|\Real  z|}\underset{\xi\in\partial\Delta (\zeta ;1)}{\sup}(1+|\xi |)^Me^{N|\Real\xi |}\\
 &\leq \tilde{C}(1+|\zeta |)^M(1+|z|)^Me^{N|\Real \zeta |}e^{N|\Real z|}.
\end{align}
\qed
\end{prf}

\begin{prop}\label{prop51}
Suppose a complex of free $\mathcal{H}$-modules, 
\begin{equation}
0\rightarrow \mathcal{H}^{1\times r_{N}}\xrightarrow{\times \mathbb{P}_N}\cdots \rightarrow \mathcal{H}^{1\times r_1}\xrightarrow{\times \mathbb{P}_1}\mathcal{H}^{1\times r_0}
\end{equation}
is given. Put $E_j={}_n\mathcal{O}^{1\times r_j}$ and consider the associated complex obtained by tensoring ${}_n\mathcal{O}$ to the above complex\footnote{Recall that $-\underset{\mathcal{H}}{\otimes }\mathcal{O}(\mathbb{C}^n)$ is an exact functor by Proposition \ref{firstprop} (3).}:
\begin{equation}\label{ThisComplex}
0\rightarrow E_N\xrightarrow{f_N}\cdots \rightarrow E_1\xrightarrow{f_1}E_0.
\end{equation}
Assume (\ref{ThisComplex}) is generically exact and equip these trivial bundles with standard Hermitian metrics and consider the associated currents $U$ and $R$. Decompose $U$ and $R$ as
\begin{equation}
U=\displaystyle \sum_{k>l}U^l_k, \;\;R=\displaystyle \sum_{k>l}R^l_k
\end{equation}
 and 
\begin{equation}
U^l_k=\displaystyle \sum_{|I|=k-l-1}U^l_{k,I}d\bar{\zeta}^I,\;\;R^l_k=\displaystyle \sum_{|I|=k-l}R^l_{k,I}d\bar{\zeta}^I,
\end{equation}
where $\zeta$ is the standard coordinate of $\mathbb{C}^n$.

Under the assumptions and notation above, there is a non-zero polynomial $\mathcal{R}(\zeta_1)$ such that $\mathcal{R}(\zeta_1)U^l_{k,I}$ and $\mathcal{R}(\zeta_1)R^l_{k,I}$ are distributions with Paley-Wiener growth, i.e., there are non-negative integers $M,\; N,\; k,$ and $k^\prime\in \mathbb{Z}_{\geq 0}$ and a positive constant $C>0$ so that for any compactly supported smooth function $\phi\in C^\infty_c(\mathbb{C}^n)$ the following inequalities hold:
\begin{equation}
\left|\langle\mathcal{R}(\zeta_1)U^l_{k,I},\phi\rangle\right|\leq C \underset{\underset{|\alpha|\leq k,|\beta |\leq k^\prime }{\zeta\in \supp\phi}}{\sup}\left|(\partial^\alpha_{\zeta }\partial^\beta_{\bar{\zeta }}\phi(\zeta ))(1+|\zeta |)^Me^{N|\Real\zeta_1 |}\right|,
\end{equation}
\begin{equation}
\left|\langle\mathcal{R}(\zeta_1)R^l_{k,I},\phi\rangle\right|\leq C \underset{\underset{|\alpha|\leq k,|\beta |\leq k^\prime }{\zeta\in \supp\phi}}{\sup}\left|(\partial^\alpha_{\zeta }\partial^\beta_{\bar{\zeta }}\phi(\zeta ))(1+|\zeta |)^Me^{N|\Real\zeta_1 |}\right|.
\end{equation}
\end{prop}

\begin{prf}
We use the same symbols as in Section 4. Let $E$ and $Q$ be trivial bundles on $\mathbb{C}^n$. We  take global holomorphic frames $\{ \epsilon_1,\ldots ,\epsilon_r\}\; (r=\rank Q)$, $\{ e_1,\ldots ,e_s\}\; (s=\rank E)$, and equip $E$ and $Q$ with hermitian metrics so that these frames are orthonormal basis at each point of $\mathbb{C}^n$. Notice that the dual section $e^*_i$ of $e_i$ is a holomorphic section of $E^*$. Take $f\in \Hom(E,Q)$ and express it as $f=\displaystyle\sum_{i,j}f_{ij}e_i^*\otimes \epsilon_j,\; f_{ij}\in\mathcal{O}(\mathbb{C}^n)$. If the optimal rank of $f$ is $\rho$, we have $F=\frac{1}{\rho !}\wedge^\rho f=\displaystyle\sum_{|I|=|J|=\rho}F_{I,J}e_I^*\otimes \epsilon_J$, where each coefficient $F_{I,J}\in\mathcal{O}(\mathbb{C}^n)$ is a product of finitely many $f_{i,j}$.
Taking into account that the standard dual $\epsilon^j$ of $\epsilon_j$ is equal to the hermitian dual $\epsilon^*_j$ of it, we have $S=\displaystyle\sum_{|I|=|J|=\rho }\bar{F}_{I,J}e_I\otimes\epsilon_J^* $. 
Since $s=(\displaystyle\sum_{i=1}^s\displaystyle\sum_{j=1}^r\delta_{f_{ij}e^*_i}\otimes\delta_{\epsilon_j})^{\rho -1}S/\rho !$, if we represent $s$ as $s=\displaystyle\sum_{i,\; j}s_{ij}e_i\otimes\epsilon_j^*$, each $s_{ij}$ is a product of finitely many $f_{ij}$ and $\bar{f_{ij}}$. Now we take $f=f_k$, $E=E_k$, and $Q=E_{k-1}$. We are going to estimate $\bar{\partial}\sigma$ for $\sigma=\frac{s}{|F|^2}$.
Writing $F_{I,J}=\frac{p_{I,J}}{\phi_{I,J}},\; \phi\in\mathbb{C}[\zeta_1],$ $p_{I,J}\in\mathbb{C}[\zeta_1,\ldots ,\zeta_n,e^{\zeta_1}]$,  letting $\phi$ be the least common multiplier of $\phi_{I,J}$ , $F_{I,J}=\frac{\tilde{F}_{I,J}}{\phi },\; \tilde{F}_{I,J}\in\mathbb{C}[\zeta_1,\ldots ,\zeta_n,e^{\zeta_1}]$ and putting  $\tilde{F}=\displaystyle\sum_{I,J}\tilde{F}_{I,J}e^*_I\otimes\epsilon_J$, we have an identity
\begin{equation}
\left(\frac{s}{|F|^{2}}\right)=|\phi |^{2}\left(\frac{s}{|\tilde{F}|^{2}}\right).
\end{equation}
%We denote the product $|\phi|^2s$ again by $s$.
By Lemma \ref{prop61}, each $f_{ij}$ satisfies the estimate 
\begin{equation}\label{Inequality}
|f_{ij}(\zeta )|\leq C(1+|\zeta |)^Me^{N|Re\zeta_1|}.
\end{equation}
Therefore we can write 
\begin{equation}
\bar{\partial}\sigma =\frac{\tilde{s}}{|\tilde{F}|^4},\; \tilde{s}=\displaystyle\sum_{i,\;j}\tilde{s}_{ij}e_i\otimes\epsilon_j^*,
\end{equation}
so that each $\tilde{s}_{ij}$ satisfies the estimate of the type (\ref{Inequality}).

In the followings, we construct the associated current $U$. Under the notation analogous to the argument above, we have a representation
\begin{equation}
u^l_{k,I}=\frac{1}{|\tilde{F}_k|^4\cdots|\tilde{F}_{l+1}|^4}\displaystyle\sum_{i,j} g_{ij}e_i^{(k)}\otimes e_j^{(l)*},
\end{equation}
where $\{e_i^{(k)}\}$ and $\{e_j^{(l)}\}$ are global holomorphic frames of $E_k$ and $E_l$ and $g_{ij}$ is a linear combination of products of holomorphic and antiholomorphic functions satisfying the estimates of the type (\ref{Inequality}). 
We rewrite the denominator as a sum
\begin{equation}
|\tilde{F}_k|^{2}\cdots|\tilde{F}_{l+1}|^2=|f_1|^2+\cdots+|f_p|^2=||f||^2,
\end{equation}
where each $f_j\in\mathbb{C}[\zeta_1,\ldots ,\zeta_n,e^{\zeta_1}]$ is a product of some entries of $\tilde{F}_{l+1},\cdots,\tilde{F}_k$. 
In view of Proposition \ref{prop42}, if $t_1,\cdots,t_p$ are positive numbers, the current $U^l_{k,I}$ is equal to
\begin{equation}
\frac{|f|^{*(t\lambda)}}{||f||^4}\displaystyle\sum g_{ij}e_i^{(k)}\otimes e_j^{(l)*}|_{\lambda=0}.
\end{equation}
%We want to find a univariate polynomial $\mathcal{R}(\zeta_1)$ so that the current $\mathcal{R}(\zeta_1)\frac{|f|^{2*(t\lambda)}}{||f||^4}\displaystyle\sum_{i,j} g_{ij}e_i^{(k)}\otimes e_j^{(l)*}|_{\lambda=0}$ is of Paley-Wiener growth. 

%If we put $|\tilde{F}|^{*(2t\lambda )}=\displaystyle\prod_{|I|=|J|=\rho }|\tilde{F}_{I,J}|^{2t_{I,J}\lambda }$, where $t_{I,J}$ are positive real numbers which will be chosen below, then we have $\left(\frac{|\tilde{F}|^{2\lambda }}{|\tilde{F}|^{2m}}\right)= |\tilde{F}|^{2\lambda }|\tilde{F}|^{*(-2t\lambda )}\left(\frac{|\tilde{F}|^{*(2t\lambda )}}{|\tilde{F}|^{2m}}\right).$ By means of toric desingularization and Lemma \ref{lem62}, one can confirm $ |\tilde{F}|^{2\lambda }|\tilde{F}|^{*(-2t\lambda )}|_{\lambda =0}=1$, hence we have the identity 
%$$g\cdot\left(\frac{|\tilde{F}|^{2\lambda }}{|\tilde{F}|^{2m}}\right)|_{\lambda =0}=g\cdot \left(\frac{|\tilde{F}|^{*(2\lambda )}}{|\tilde{F}|^{2m}}\right)|_{\lambda =0}.$$
A priori, this current is holomorphic at $\lambda=0$. Therefore, the resulting current is a combination of the constant term of the Taylor expansion of $\frac{|f|^{*(t\lambda)}}{||f||^4}$ multiplied by some products of holomorphic and anti-holomorphic functions of Paley-Wiener growth. 

Now we need the following lemma.

\end{prf}

\begin{lem}
Let $T\in\mathcal{D}^\prime (\mathbb{C}^n)$ be a distribution with Paley-Wiener growth. If $g\in\mathcal{O}(\mathbb{C}^n)$ is a holomorphic function which satisfies the estimate
\begin{equation}\label{PWEstimate}
|g(\zeta )|\leq C(1+|\zeta |)^Me^{N|\Real\zeta_1|}
\end{equation}
for any $\zeta\in\mathbb{C}^n,$ then two distributions $gT$ and $\bar{g}T\in\mathcal{D}^\prime (\mathbb{C}^n)$ both have Paley-Wiener growth. 
\end{lem}

\begin{prf}
We prove the lemma for $gT$ since the argument for $\bar{g}T$ is similar. Take any $\phi\in\mathcal{C}^\infty_c(\mathbb{C}^n)$. By definition, 
\begin{equation}
|\langle gT,\phi \rangle|=|\langle T,g\phi \rangle|\leq C\underset{\underset{|\alpha|\leq k,|\beta |\leq k^\prime }{\zeta\in supp\phi}}{\sup}\left|(\partial^\alpha_{\zeta }(g(\zeta )\partial^\beta_{\bar{\zeta }}\phi(\zeta )))(1+|\zeta |)^Me^{N|\Real\zeta_1|}\right|.
\end{equation}
It can be verified from Cauchy's integral representation theorem that derivatives of $g$ satisfy the estimates of the type (\ref{PWEstimate}). Combining this observation with the inequality above, we have the lemma.\qed
\end{prf}

\noindent
{\bf End of proof of Proposition \ref{prop51}}
Thanks to this lemma, the problem is reduced to showing that the constant term of the Taylor expansion of $\frac{|f|^{*(t\lambda)}}{|f|^4}$ around the origin is of Paley-Wiener growth. The Theorem \ref{Berenstein} below  says that the constant term is indeed of Paley-Wiener growth if it is multiplied by a non-zero polynomial $\mathcal{R}(\zeta_1)$. The proof is completed for $U$. 
For $R$, we can observe that an identity  $\bar{\partial}U^l_k=R^l_k+(|f|^{*(t\lambda)}\bar{\partial }u^l_k)|_{\lambda =0}$ holds. 
Using this identity, one can show that there is a univariate polynomial $\mathcal{R}(\zeta_1)$ so that $\mathcal{R}(\zeta_1)R^l_k$ is of Paley-Wiener growth in the same manner as we proved it for $U$.\qed
%Since left hand side is of Paley-Wiener growth and we can check that the second term of the right hand side is of Paley-Wiener growth as we confirmed $U^l_k$ is of Paley-Wiener growth, we can conclude $R^l_k$ is of Paley-Wiener growth.  \qed

\begin{thm}[{\cite[Proposition 3.2]{BY1}}]\label{Berenstein}
Let $f_1,\cdots,f_p\in\mathbb{C}[\zeta_1,\cdots,\zeta_n,e^{\zeta_1}]$, where p is a positive integer, then, for any $t\in(0,1)^p$ outside a countable union of algebraic hypersurfaces and any $m\in\mathbb{Z}_{>0}$, there is a polynomial $\mathcal{R}(\zeta_1)$ and constants $C>0,\; M,\; N,\; k,\; k^\prime\in\mathbb{Z}_{>0}$ such that if $a_j\in\mathcal{D}^\prime (\mathbb{C}^n)$ denote the coefficients of the Taylor expansion
\begin{equation}
\frac{|f|^{*(t\lambda)}}{\| f\|^{2m}}=\displaystyle\sum_{j=0}^\infty a_j\lambda^j,
\end{equation}
where $|f|^{*(t\lambda)}=|f_1|^{t_1\lambda}\cdots |f_p|^{t_p\lambda}$ and $\| f \| ^2=|f_1|^2+\cdots +|f_p|^2$, then for  $\phi\in C^\infty_c(\mathbb{C}^n)$,
\begin{equation}
|\langle\mathcal{R}(\zeta_1)a_0,\phi\rangle|\leq C\sup_{\underset{|\alpha|\leq k, |\beta|\leq k^\prime}{\zeta\in\supp\phi}}\left|\left( \partial^\alpha_\zeta\partial^\beta_{\bar{\zeta}}\phi(\zeta)\right)(1+|\zeta|)^M e^{N|\Real\zeta_1|}\right|.
\end{equation}
\end{thm}

\begin{ex}
We consider a univariate function $e^\zeta -1\in\mathcal{O}(\mathbb{C}_\zeta)$ and the associated generically exact complex of trivial vector bundles 
\begin{equation}\label{SimpleComplex}
0\rightarrow\mathcal{O}\overset{e^\zeta-1}{\rightarrow}\mathcal{O}\rightarrow 0.
\end{equation}
By definition, the current $U$ associated to the complex (\ref{SimpleComplex}) is the principal value current ${\rm v. p.}\;\frac{1}{e^\zeta-1}$. A direct computation leads to the following formula:
\begin{align}
{\rm v. p.}\;\frac{1}{e^\zeta-1}&= e^{-\zeta}\frac{\partial}{\partial\zeta}\left(\log|e^\zeta-1|^2\right)\\
 &=e^{-\zeta}\frac{\partial}{\partial\zeta}\left(e^{-\zeta}\frac{\partial}{\partial\zeta}\left((e^\zeta-1)\log|e^\zeta-1|^2\right)\right).
\end{align}
Since $(e^\zeta-1)\log|e^\zeta-1|^2$ is locally bounded and its growth order at infinity is exponential 1, we can confirm that ${\rm v. p.}\;\frac{1}{e^\zeta-1}$ is of Paley-Wiener growth. We can also compute the residue current $R$ explicitly:
\begin{equation}
R=\pi\displaystyle\sum_{m\in\mathbb{Z}}\delta(\zeta-2\pi\sqrt{-1}m)d\bar{\zeta}.
\end{equation}
This formula also tells us that $R$ is indeed of Paley-Wiener growth.
\end{ex}

Now we can finally prove the division with bounds. The essential point is the construction of the suitable weight in view of Theorem \ref{intrep}. We first work on the growth condition which arises from Fourier-Borel transform of analytic functionals.

For any compact convex subset $K$ of $\mathbb{C}^n$ with smooth boundary, and $h\in\mathcal{O}(\mathbb{C}^n)^{1\times r_0}$, we put 
\begin{equation}
|h|_{K}=\underset{\zeta\in\mathbb{C}^n}{\sup}e^{-H_K(\zeta )}|h(\zeta )|, \;H_K(\zeta )=\underset{z\in K}{\sup}\Real\langle z,\zeta \rangle,\;\zeta\in\mathbb{C}^n,
\end{equation}
where $\langle\cdot ,\cdot \rangle$ is the standard Euclidian product.
$H_K$ is smooth except at the origin so we smoothen out $H_K(\zeta )$ around $\zeta =0$ so that the resulting function is convex. This function is denoted by $\rho_K$.
We put 
\begin{equation}
\rho_K^\prime (\zeta )=\left(2\frac{\partial\rho_K}{\partial\zeta_1}(\zeta ),\ldots ,2\frac{\partial\rho_K}{\partial\zeta_n}(\zeta )\right),
\end{equation}
and 
\begin{equation}
g_K (\zeta ,z)=\exp\{ -\langle\rho_K^\prime (\zeta),\zeta -z\rangle-\frac{\sqrt{-1}}{\pi }\partial\bar{\partial }\rho_K\}=\exp\{ \nabla_{\zeta -z}(\frac{-\partial\rho_K}{\pi\sqrt{-1}})\}.
\end{equation}
It can easily be seen that $(g_K)_{0,0}(z,z)=1$ and $\nabla_{\zeta -z}g_K =0$. Furthermore, since $\rho_K$ is convex, we have an inequality
\begin{equation}
|\exp\{-\langle\rho_K^\prime (\zeta),\zeta -z\rangle\}|\leq \exp\{ \rho_K (z)-\rho_K(\zeta)\}\;\; (\zeta,  z\in\mathbb{C}^n).
\end{equation}

In the same manner, we put $\tilde{\rho} (\zeta )\overset{def}{=}\tilde{\rho}(\xi_1 )\overset{def}{=}|\xi_1|,$ where $\zeta =\xi +\sqrt{-1}\eta $.
We smoothen out $\tilde{\rho}$ around $\xi_1 =0$ so that the resulting function is  convex and this function is still denoted by $\tilde{\rho}$.
Under the similar notation as above, we have identities $(\tilde{g})_{0,0}(z,z)=1$, $\nabla_{\zeta -z}\tilde{g} =0,$
and an inequality
\begin{equation}
|\exp\{ -\langle\tilde{\rho}^\prime (\zeta),\zeta -z\rangle\}|\leq \exp\{ \tilde{\rho} (z)-\tilde{\rho} (\zeta)\},\;\;(\zeta, z\in\mathbb{C}^n).
\end{equation}
%We put $\tilde{\rho}^\prime (\zeta )=\tilde{\rho}^\prime (\xi_1 )=(\frac{\partial\tilde{\rho}}{\partial\zeta_1}(\zeta ),0,\ldots ,0)$ and $\tilde{g} (\zeta ,z)=\exp\{ \langle\tilde{\rho}^\prime (\zeta),\zeta -z\rangle+\frac{\sqrt{-1}}{\pi }\partial\bar{\partial }\tilde{\rho} \}=\exp\{ \nabla_{\zeta -z}(\frac{\partial\tilde{\rho}}{\pi\sqrt{-1}})\}$.
%Then we have $(\tilde{g})_{0,0}(z,z)=1,\;\bigtriangledown_{\zeta -z}\tilde{g} =0$.\\
%Furthermore, $|\exp\{ \langle\tilde{\rho}^\prime (\zeta),\zeta -z\rangle\}|\leq \exp\{ \tilde{\rho} (z)-\tilde{\rho} (\zeta)\}\;\;\forall \zeta, \forall z\in\mathbb{C}^n$

%\begin{thm}
%$\Omega\subset\mathbb{R}^n$; open convex subset such that $\mathbb{Z}\mathbb{e}_1+\Omega=\Omega$ and $\partial\Omega$ is smooth. We put $\mathcal{O}(\mathbb{C}^n)_p^{1\times r_0}=\{ h\in\mathcal{O}(\mathbb{C}^n)^{1\times r_0}|\exists M\in\mathbb{Z}_{\geq 0},\;\exists K\subset\Omega ;compact\; convex,\; |h|_{M,K}<\infty \}$ For any $\mathbb{P}\in M(r_1,r_0;\mathcal{H})$ we have\\
%$\mathcal{O}(\mathbb{C}^n)_p^{1\times r_0}\cap (\mathcal{O}(\mathbb{C}^n)^{1\times r_1}\cdot\mathbb{P})=\mathcal{O}(\mathbb{C}^n)_p^{1\times r_1}\cdot\mathbb{P}$
%\end{thm}

We have now reached one of principal results of this paper: the division with bounds. We denote by $\mathcal{O}(\mathbb{C}^n)_p^{1\times r_0}$ the set of $r_0$ vectors $h$ with entries in $\mathcal{O}(\mathbb{C}^n)$ such that $|h|_K<\infty$ for some compact convex subset $K\subset\Omega$.

\begin{thm}\label{cpxdivision}
Let $\Omega$ be an open convex subset of $\mathbb{C}^n$ such that $\mathbb{Z}\mathbf{e}_1+\Omega=\Omega$. Then for any matrix $\mathbb{P}\in M(r_1,r_0;\mathcal{H})$ of generalized $D\Delta$-operators, one has the fundamental identity:
\begin{equation}
\mathcal{O}(\mathbb{C}^n)_p^{1\times r_0}\cap (\mathcal{O}(\mathbb{C}^n)^{1\times r_1}\cdot\mathbb{P})=\mathcal{O}(\mathbb{C}^n)_p^{1\times r_1}\cdot\mathbb{P}.
\end{equation}
\end{thm}

For the proof, we need a simple algebraic lemma which will reduce the problem from the division with bounds for a submodule to that for an ideal.

\begin{lem}
If for any finitely generated ideal $I$ of $\mathcal{H}$, one has the identity 
\begin{equation}
\left(\mathcal{O}(\mathbb{C}^n)\cdot I\right)\cap\mathcal{O}(\mathbb{C}^n)_p=\mathcal{O}(\mathbb{C}^n)_p\cdot I,
\end{equation}
then for any positive integer $r\in\mathbb{Z}_{>0}$, and for any finitely generated submodule $N$ of $\mathcal{H}^{1\times  r}$, one has the identity 
\begin{equation}
\left(\mathcal{O}(\mathbb{C}^n)\cdot N\right)\cap\mathcal{O}(\mathbb{C}^n)_p^{1\times r}=\mathcal{O}(\mathbb{C}^n)_p\cdot N.
\end{equation}
\end{lem}

\begin{prf}
First of all, notice that it is equivalent to proving the following claim to prove this lemma:\\
(Claim) For any coherent $\mathcal{H} $ module $M$, the complex
\begin{equation}
0\rightarrow\mathcal{O}(\mathbb{C}^n)_p\otimes_\mathcal{H}M\rightarrow\mathcal{O}(\mathbb{C}^n)\otimes_\mathcal{H}M
\end{equation}
is exact.

\noindent
In fact, if we put $M=\mathcal{H}^{1\times r}/N$ in the claim, we have an inclusion $\mathcal{O}(\mathbb{C}^n)\cdot N\cap\mathcal{O}(\mathbb{C}^n)_p^{1\times r}\subset\mathcal{O}(\mathbb{C}^n)_p\cdot N$.
Since the other inclusion is obvious, we get the lemma.

Now we prove the claim by the induction on the number of generators of $M$.
When $M$ is cyclic, we can assume there is a finitely generated ideal $I$ of $\mathcal{H}$ such that $M\simeq\mathcal{H}/I$ so the claim follows by the assumption of lemma. Assume that the claim is valid when the number of generators of $M$ is less than $n$. Let $M$ be a coherent $\mathcal{H}$-module generated by $n$ elements. We can find an exact sequence $0\rightarrow M^\prime\rightarrow M\rightarrow M/M^\prime\rightarrow 0$, where $M^\prime$ is generated by at most $n-1$ elements and $M/M^\prime$ is cyclic. Consider the following diagram.
\[
\xymatrix{
 & 0 \ar[d] &  & 0 \ar[d]  & \\
 & \mathcal{O}(\mathbb{C}^n)_p\underset{\mathcal{H}}{\otimes} M^\prime \ar[r] \ar[d] & \mathcal{O}(\mathbb{C}^n)_p\underset{\mathcal{H}}{\otimes} M \ar[r] \ar[d] & \mathcal{O}(\mathbb{C}^n)_p\underset{\mathcal{H}}{\otimes} M/M^\prime \ar[r] \ar[d] & 0\\
0 \ar[r] & \mathcal{O}(\mathbb{C}^n)\underset{\mathcal{H}}{\otimes} M^\prime \ar[r]  & \mathcal{O}(\mathbb{C}^n)\underset{\mathcal{H}}{\otimes} M \ar[r]  & \mathcal{O}(\mathbb{C}^n)\underset{\mathcal{H}}{\otimes} M/M^\prime \ar[r]  & 0.\\
}
\]
The first and the second row are exact by the exactness of tensor and Proposition \ref{firstprop} (3), and the first and the third column are exact by the inductive assumption. We can now apply the snake lemma to conclude that the middle vertical map is injective.\qed
\end{prf}

\noindent
{\bf Proof of Theorem \ref{cpxdivision}}
We basically follow the argument of \cite{BY1}. The crucial difference is that we make use of a more general division formula than in \cite{BY1}. By the previous lemma, we can assume $r_0=1$ and we write $\mathbb{P}=\;{}^t(P_1,\cdots ,P_{r_1})$.
Let $h\in\mathcal{O}(\mathbb{C}^n)_p\cap (\mathcal{O}(\mathbb{C}^n)^{1\times r_1}\cdot\mathbb{P})$ such that $|h|_{K}<\infty $ for a compact convex set $K$. We can assume $\partial K$ is smooth since there is always a smooth convex function $\varphi$ on $\Omega$ such that $\{\varphi<a\}$ is relatively compact in $\Omega$ and $K\subset\{\varphi<0\}$. The existence of such smooth convex function is verified in the same manner as that of plurisubharmonic function $\psi$ which exhausts a given pseudoconvex domain $D$ and satisfies that $L\subset\{\psi<0\}$ for given holomorphically convex compact subset $L$ of $D$. See Theorem 5.1.6 in \cite{Ho}.
By Theorem \ref{thm53} we can take a finite free resolution of finite length
\begin{equation}
0\rightarrow  \mathcal{H}^{1\times r_{N}}\xrightarrow{\times \mathbb{P}_N}\cdots \rightarrow \mathcal{H}^{1\times r_1}\xrightarrow{\times \mathbb{P}_1=\mathbb{P}}\mathcal{H}\;(exact).
\end{equation}
We can also take Hefer forms for this exact sequence as in Proposition \ref{prop55}, and choose a polynomial $\mathcal{R}(z_1)$ as in Proposition \ref{prop51}. We denote its distinct roots by $\alpha_1,\cdots ,\alpha_k$, and let $\mu_1,\cdots ,\mu_k$ be their respective multiplicities.

We can expand each $P_j$ as a power series of $z_1-\alpha_l$ whose coefficients are polynomials in $z^\prime =(z_2,\cdots ,z_n)$. If we truncate this series at the term $(z_1-\alpha_l)^{\mu_l}$, we obtain a polynomial $P_{j,l}$. 
By the classical Ehrenpreis fundamental principle (or division with bounds, see \cite[Theorem 4.2]{Eh} or \cite[IV, \S 5, $5^\circ$ The fundamental theorem]{P}),  possibly replacing $K$ by its compact convex neighbourhood, we have a representation
\begin{equation}
h=\displaystyle\sum_{j=1}^{r_1}G_{j,l}P_{j,l}+(z_1-\alpha_l)^{\mu_l}G_{r_1+1,l},\;\; |G_{j,l}|_K<\infty .
\end{equation}
Now for each fixed $z^\prime$ we can construct a polynomial $G_{j}(z_1,z^\prime )$ in $z_1$, by means of Lagrange interpolation formula such that for each $l$,
\begin{equation}
G_{j}(z_1,z^\prime )-G_{j,l}(z_1,z^\prime)=O((z_1-\alpha_l)^{\mu_l}).
\end{equation}
The explicit interpolation formula implies that $|G_j|_K<\infty$.

Now by means of these functions $G_j$, we have 
\begin{equation}
h-\displaystyle\sum_{j}^{r_1}G_jP_j=h-\displaystyle\sum_{j}^{r_1}G_{j,l}P_j+\displaystyle\sum_{j}^{r_1}(G_{j,l}-G_j)P_j=O((z_1-\alpha_l)^{\mu_l}).
\end{equation}
We can conclude that there is a representation 
\begin{equation}
h=\displaystyle\sum_{j=1}^{r_1}G_jP_j+\mathcal{R}(z_1)G_{r_1+1}.
\end{equation}
Note here that we have an inequality $|\mathcal{R}(z_1)G_{r_1+1}|_K<\infty$ again by possibly replacing $K$ by its compact convex neighbourhood. This implies $|G_{r_1+1}|_K<\infty$ in view of P\'olya-Ehrenpreis-Malgrange division lemma. (This is nothing but the division with bounds for principal ideal.)
We can confirm that $\mathcal{R}(z_1)G_{r_1+1}$ is in the ideal generated by $P_1,\cdots ,P_{r_1}$ so we want to give bounds for coefficients of each $P_j$ by means of the residue current.

Choose a function $\chi\in C^\infty (\mathbb{R},\mathbb{R})$  so that
\begin{equation}
\chi (t)= 1\; (\text{if }|t|<1), \;\;\chi (t)=0\; (\text{if } |t|>2).
\end{equation}
Putting
\begin{equation}
\chi_k (\zeta )=\chi (\frac{|\zeta |}{k}),
\end{equation}
we introduce a weight 
\begin{equation}
g_k=\chi_k -\partial\chi_k\wedge \left(\frac{b}{\nabla_{\zeta -z}b}\right),
\end{equation}
where $b$ is a (1,0)-form defined by the formula
\begin{equation}
b=\frac{1}{2\pi\sqrt{-1}}\frac{\partial |\zeta-z|^2}{|\zeta-z|^2}.
\end{equation}
Here, $\frac{b}{\nabla_{\zeta -z}b}$ is M. Andersson's notation in \cite{AII}, which reads
\begin{equation}
\frac{b}{\nabla_{\zeta -z}b}=b+b\wedge\bar{\partial}b+\dots +b\wedge(\bar{\partial}b)^{n-1}.
\end{equation}
Note that $g_k$ is a smooth weight when $|z|<k$. We further put 
\begin{equation}
g(\zeta ,z)=\frac{1+\langle\bar{\zeta},z\rangle}{1+|\zeta |^2}+\frac{\sqrt{-1}}{2\pi}\partial\bar{\partial}\log(1+|\zeta |^2).
\end{equation}
It can be verified that $g$ satisfies $g_{0,0}(z,z)=1$ and $\nabla_{\zeta -z}g=0$.

Now consider a current $g_k\wedge g_K\wedge g^\mu\wedge\tilde{g}^\nu$ for non-negative integers $\mu$ and $\nu$.
This is a weight with respect to $z$ when $|z|<k$. By Proposition \ref{prop46}, we have
\begin{equation}
\mathcal{R}(z_1)G_{r_1+1}(z)=\left(\int_\zeta H^1U\mathcal{R}(\zeta_1)G_{r_1+1}(\zeta)\wedge g_k\wedge g_K\wedge g^\mu\wedge \tilde{g}^\nu\right)\cdot\mathbb{P}(z),\; |z|\leq k.
\end{equation}
Let us observe that $g=o((1+|\zeta|)^{-1})$ for each fixed $z$.
If we take suitably large $\mu$, $\nu$, since $\mathcal{R}(z_1)U$ is of Paley-Wiener growth, each derivative of $g$ and $\tilde{g}$ satisfies an estimate of Paley-Wiener type, and $g_k\rightarrow 1$ and since  $\bar{\partial}g_k\rightarrow 0$ in $C^\infty (\mathbb{C}^n)$ as $k\rightarrow\infty$, we can observe that the integral 
\begin{equation}
\int_\zeta H^1U\mathcal{R}(\zeta_1)G_{r_1+1}(\zeta)\wedge g_K\wedge g^\mu\wedge \tilde{g}^\nu
\end{equation}
converges so that we have
\begin{equation}
\mathcal{R}(z_1)G_{r_1+1}(z)=\left(\int_\zeta H^1U\mathcal{R}(\zeta_1)G_{r_1+1}(\zeta)\wedge g_K\wedge g^\mu\wedge\tilde{g}^\nu\right)\cdot\mathbb{P}(z),
\end{equation}
and 
\begin{equation}
\left|\int_\zeta H^1U\mathcal{R}(\zeta_1)G_{r_1+1}(\zeta)\wedge g_K\wedge  g^\mu\wedge\tilde{g}^\nu \right|\leq C(1+|z|)^{M^\prime}e^{\rho (z)}e^{\nu\tilde{\rho}(z)}.
\end{equation}
Since for the convex hull $\tilde{K}={\rm c.h.}(K\pm\nu\mathbf{e}_1)\subset\Omega$, we have $\rho (z)+\nu\tilde{\rho} (z)\leq H_{\tilde{K}}(z)$ as long as $z$ is outside a small neighbourhood of the origin and $(1+|z|)^{M^\prime}\leq\tilde{C}e^{\varepsilon |z|}$ for any $\epsilon>0$, we can conclude that  $h\in\mathcal{O}(\mathbb{C}^n)_p^{1\times r_1}\cdot\mathbb{P}$ holds. \qed

\vspace{1em}

In  the same manner, we can solve the division with bounds where growth conditions arise from the Fourier transform. For a convex compact subset $K$ of $\mathbb{R}^n$, a non-negative integer $M\in\mathbb{Z}_{\geq 0}$, and a vector of holomorphic functions $h\in\mathcal{O}(\mathbb{C}^n)^{1\times r_0}$, we put 
\begin{equation}
|h|_{M,K}=\underset{\zeta\in\mathbb{C}^n}{\sup}(1+|\zeta |)^{-M}e^{-H_K(\zeta )}|h(\zeta )|\end{equation}
 and 
\begin{equation}
H_K(\zeta )=\underset{x\in K}{\sup}\langle x,\eta \rangle,\;\zeta =\xi +\sqrt{-1}\eta\in\mathbb{C}^n,
\end{equation}
where $\langle\cdot ,\cdot \rangle$ is the standard Euclidian product.\\

We use a different embedding $\mathcal{H}\subset\mathcal{O}(\mathbb{C}^n_\zeta )$ from that in Section 3 by putting $z_i\mapsto -\sqrt{-1}\zeta_i $, $e^{z_1}\mapsto e^{-\sqrt{-1}\zeta_1}$. Note that under this embedding, all the results in Section 3 still hold. We denote by $\mathcal{O}(\mathbb{C}^n)_{p^\prime}^{1\times r_0}$ the set of vectors $h$ with entries in $\mathcal{O}(\mathbb{C}^n)$ such that $|h|_{M,K}<\infty$ for some positive integer $M$ and a compact convex subset $K\subset\Omega$.

\begin{thm}\label{division}
Let $\Omega$ be an open convex subset of $\mathbb{R}^n$ such that $\mathbb{Z}\mathbf{e}_1+\Omega=\Omega$. Then for any matrix $\mathbb{P}\in M(r_1,r_0;\mathcal{H})$ of generalized $D\Delta$-operators, one  always has the identity
\begin{equation}
\mathcal{O}(\mathbb{C}^n)_{p^\prime}^{1\times r_0}\cap (\mathcal{O}(\mathbb{C}^n)^{1\times r_1}\cdot\mathbb{P})=\mathcal{O}(\mathbb{C}^n)_{p^\prime}^{1\times r_1}\cdot\mathbb{P}.
\end{equation}
\end{thm}

\begin{prf}
Since the proof is parallel to that of Theorem \ref{cpxdivision}, we only give definitions of the weights. Let us take any compact convex subset $K$ of $\Omega$ and define $\rho_K(\zeta)=\rho(\eta)$ to be the smoothened version of the convex support function $H_K(\zeta)$. We put
\begin{equation}
g_K(\zeta,z)=\exp\{\nabla_{\zeta-z}(-\frac{\partial\rho_K}{\pi\sqrt{-1}})\}.
\end{equation}
Similarly, we put $\tilde{\rho}(\zeta)\overset{def}{=}\tilde{\rho}(\eta_1)\overset{def}{=}|\eta_1|$ and again smoothen this out around the origin. Defining $\tilde{g}$ by $\tilde{g}=\exp\{\nabla_{\zeta-z}(-\frac{\partial\tilde{\rho}}{\pi\sqrt{-1}})\}$, we can mimic the proof of Theorem \ref{cpxdivision}.\qed
\end{prf}

\section{Ehrenpreis-Malgrange type theorems for generalized $D\Delta$ equations}

%We here give an elementary version of injectivity result for exponential polynomial functions.
%\begin{prop}\ 
%We put $E=\{f(z)|f(z)\; is\; a\; linear\; combination\; of\; p(z)e^{\alpha \cdot z},\; z\in\mathbb{C}^n, \alpha\in\mathbb{C}^n\}$ where $\cdot$ stands for the dot product.\\
%Then for any coherent $\mathcal{H}$ module $M$, $Ext_\mathcal{H}^i(M,E)=0\; \forall i>0$
%\end{prop}

%We finally investigate the system of differential equations with constant coefficients and commensurate time lags. \\

We can finally establish Ehrenpreis-Malgrange type theorems for $D\Delta$-equations.
Firstly, we give the simplest version of our main theorem which does not require any topological consideration.

Let us begin with defining the action of $\mathcal{H}$ on holomorphic functions.
Let $\Omega\subset\mathbb{C}^n$ be an open set which satisfies $\Omega +\mathbb{R}\mathbf{e}_1=\Omega $, that is, an open set which is an inverse image of another open set with respect to the natural projection $\pi :\mathbb{C}^n\rightarrow Y=\mathbb{C}^n/\mathbb{R}\mathbf{e}_1\simeq \sqrt{-1}\mathbb{R}\times\mathbb{C}^{n-1}$.
Each connected component $\Omega_i$ of $\Omega $ again satisfies $\Omega_i +\mathbb{R}\mathbf{e}_1=\Omega_i $. On each $\Omega_i$, we can define the action of $\mathcal{H}$ on $\mathcal{O} (\Omega_i )$ as in Section \ref{1dim}.

Let $\mathbf{E}$ be the subring of $\mathcal{O}(\mathbb{C}^n)$ generated by $\mathbb{C}[z]$ and $e^{\alpha z}$ ($\alpha\in\mathbb{C}$). We can decompose $\mathbf{E}$ as a $\mathbb{C}$-vector space:
\begin{equation}
\mathbf{E}=\displaystyle\bigoplus_{\alpha\in\mathbb{C}}\mathbb{C}[z]e^{\alpha z}.
\end{equation}
We put $\mathbf{E}_\alpha=\mathbb{C}[z]e^{\alpha z}$.
We can easily see that $\mathbf{E}_\alpha$ is an $\mathcal{H}$ submodule of $\mathcal{O}(\mathbb{C}^n)$.

\begin{prop}
For any positive integer $i$ and any coherent $\mathcal{H}$-module $M$ one has the following vanishing result:
\begin{equation}
\Ext^i_\mathcal{H}(M,\mathbf{E})=0
\end{equation}
\end{prop}

\begin{prf}
By the coherency of $\mathcal{H}$, it is enough to prove that for any short exact sequence
\begin{equation}
\mathcal{H}^{1\times r}\xrightarrow{\times\mathbb{Q}(z ,\sigma )}\mathcal{H}^{1\times s}\xrightarrow{\times\mathbb{P}(z ,\sigma )}\mathcal{H}^{1\times t},
\end{equation}
the associated complex obtained by applying $Hom_\mathcal{H}(-,\mathbf{E})$
\begin{equation}\label{ASequence}
\mathbf{E}^{t\times 1}\xrightarrow{\mathbb{P}(\partial ,\sigma)\cdot }\mathbf{E}^{s\times 1}\xrightarrow{\mathbb{Q}(\partial ,\sigma)\cdot}\mathbf{E}^{r\times 1}
\end{equation}
 is exact. Furthermore, it is enough to prove the statement above replacing $\mathbf{E}$ by $\mathbf{E}_0=\mathbb{C}[z]$.

We consider a perfect pairing
\begin{equation}
\langle\cdot ,\cdot\rangle :\mathbb{C}[z]\times\mathbb{C}[[z]]\rightarrow\mathbb{C}
\end{equation}
defined by
\begin{equation}
\langle  p,f\rangle =\displaystyle\sum_{n=0}^\infty n!p_nf_n,\;\; p=\displaystyle\sum_{n}p_nz^n\in\mathbb{C}[z],\; f=\displaystyle\sum_{n}f_nz^n\in\mathbb{C}[[z]].
\end{equation}
One can easily show that for any $p\in\mathbb{C}[z]$, any  $f\in\mathbb{C}[[z]]$, and any positive integer $i$,
\begin{equation}
\langle\frac{\partial}{\partial z_i}p,f\rangle =\langle p,z_if\rangle\;\;\text{and}\;\;\langle\sigma p,f\rangle =\langle p,e^{z_1}f\rangle.
\end{equation}
Therefore, the exactness of (\ref{ASequence}) is equivalent to that of
\begin{equation}
\mathbb{C}[[z]]^{1\times r}\xrightarrow{\times\mathbb{Q}(z ,e^{z_1} )}\mathbb{C}[[z]]^{1\times s}\xrightarrow{\times\mathbb{P}(z ,e^{z_1} )}\mathbb{C}[[z]]^{1\times t}.
\end{equation}
The last complex is actually exact since $\mathcal{H}\subset\mathcal{O}(\mathbb{C}^n)\subset\mathbb{C}[[z]]$ is a tower of flat extensions by Proposition \ref{firstprop} (3) and \cite[Result 6.1.7]{OB}.\qed
\end{prf}

Next, we proceed to our main result. We use some basic terminologies of topological vector spaces. See e.g., \cite{K3} or \cite[V, \S 1]{P}. First of all, remember that we can define the Fourier transform of any distribution with compact support $T\in\mathcal{E}^\prime (\mathbb{R}^n)$ by the formula $(\mathcal{F}T)(\xi )=\langle T(x),e^{-\sqrt{-1}\xi x}\rangle$.
In the same manner, for any analytic functional $T\in\mathcal{O}^\prime(\mathbb{C}^n)$, we can define the Fourier-Borel transform of $T$ by the formula $(\mathcal{F}^\mathcal{B}T)(\zeta )=\langle T(z),e^{\zeta z}\rangle$. The classical Paley-Wiener-Schwartz theorem and Ehrenpreis-Martineau theorem state that these are linear topological isomorphisms.

\begin{thm}[{\cite[V, \S 3, Proposition 2]{P}}, {\cite[Theorem 6.4.5]{Mor}}]
For any open convex subset $\Omega$ of $\mathbb{R}^n$, the Fourier transform $\mathcal{F}$ gives rise to a linear topological isomorphism
\begin{equation}
\mathcal{E}^\prime (\Omega)\overset{\mathcal{F}}{\widetilde{\rightarrow}}\mathcal{O}(\mathbb{C}^n)_{p^\prime},
\end{equation}
where $\mathcal{O}(\mathbb{C}^n)_{p^\prime}$ is equipped with the standard (DFS) topology.

Similarly, for any convex open subset $\Omega$ of $\mathbb{C}^n$, the Fourier-Borel transform gives rise to a linear topological isomorphism\\
\begin{equation}
\mathcal{O}^\prime (\Omega)\overset{\mathcal{F}^\mathcal{B}}{\widetilde{\rightarrow}}\mathcal{O}(\mathbb{C}^n)_p,
\end{equation}
where $\mathcal{O}(\mathbb{C}^n)_p$ is equipped with the standard (DFS) topology.
\end{thm}

In the following argument, we assume that $\Omega\subset\mathbb{C}^n$ is convex. Consider any short exact sequence
\begin{equation}
\mathcal{H}^{1\times r}\xrightarrow{\times\mathbb{Q}(z ,\sigma )}\mathcal{H}^{1\times s}\xrightarrow{\times\mathbb{P}(z ,\sigma )}\mathcal{H}^{1\times t}.
\end{equation}
Since $-\otimes_{\mathcal{H}}\mathcal{O}(\mathbb{C}^n)$ is exact by Proposition \ref{firstprop} (3) , we have an exact sequence
\begin{equation}
\mathcal{O}(\mathbb{C}^n)^{1\times r}\xrightarrow{\times\mathbb{Q}(\zeta,e^{\zeta_1} )}\mathcal{O}(\mathbb{C}^n)^{1\times s}\xrightarrow{\times\mathbb{P}(\zeta,e^{\zeta_1} )}\mathcal{O}(\mathbb{C}^n)^{1\times t}.
\end{equation}
In view of Fourier-Borel transform, one can show that the complex
\begin{equation}
\mathcal{O}(\mathbb{C}^n)_p^{1\times r}\xrightarrow{\times\mathbb{Q}(\zeta,e^{\zeta_1} )}\mathcal{O}(\mathbb{C}^n)_p^{1\times s}\xrightarrow{\times\mathbb{P}(\zeta,e^{\zeta_1} )}\mathcal{O}(\mathbb{C}^n)_p^{1\times t}
\end{equation}
be exact and $\mathbb{P}(\zeta ,e^{\zeta_1})$ have a closed image implies the exactness of \\
\begin{equation}
\mathcal{O} (\Omega)^{t\times 1}\xrightarrow{\mathbb{P}(\partial ,\sigma)\cdot }\mathcal{O} (\Omega)^{s\times 1}\xrightarrow{\mathbb{Q}(\partial ,\sigma)\cdot}\mathcal{O} (\Omega)^{r\times 1}.
\end{equation}
The proof of the claim above can be found in, for example, \cite[Theorem VII.1.3]{K2}. See also \cite[V, \S 1, Proposition 8]{P}.  The exactness of the former complex is satisfied by Theorem \ref{cpxdivision}.
To prove the closed range property, take any sequence $\{\mathbf{f_j}\}_j\subset\mathcal{O}(\mathbb{C}^n)_p^{1\times s}$ such that $\{\mathbf{f_j}\mathbb{P}(\zeta,e^{\zeta_1})\}_j\subset\mathcal{O}(\mathbb{C}^n)_p^{1\times t}$ converges. Then, we can see that the limit belongs to $\mathcal{O}(\mathbb{C}^n)^{1\times s}\mathbb{P}(\zeta,e^{\zeta_1})\cap\mathcal{O}(\mathbb{C}^n)_p^{1\times t}$ since  $\mathcal{O}(\mathbb{C}^n)^{1\times s}\mathbb{P}(\zeta,e^{\zeta_1})\subset\mathcal{O}(\mathbb{C}^n)^{1\times t}$ is closed. See \cite[Result 6.1.8]{OB}. Thus the limit belongs to  $\mathcal{O}(\mathbb{C}^n)_p^{1\times s}\mathbb{P}(\zeta,e^{\zeta_1})$ by Theorem\ref{cpxdivision}.

Since $\mathcal{H}$ is a coherent ring, we have proved

\begin{thm}\label{thm71}
For any open convex subset $\Omega$ of $\mathbb{C}^n$ such that $\mathbb{R}\mathbf{e}_1 +\Omega =\Omega $, for any coherent $\mathcal{H}$-module $M$, and for any positive integer $i$, one has the following vanishing result:
\begin{equation}
\Ext^i_\mathcal{H}(M,\mathcal{O} (\Omega ))=0.
\end{equation}
\end{thm}

We can define the action of $\mathcal{H}$ on $C^\infty$ functions in the same way as we defined one on $\mathcal{O}$, and we can prove the following theorem in view of Fourier transform.

\begin{thm}\label{thm72}
For any open convex subset $\Omega$ of $\mathbb{R}^n$ such that $\mathbb{R}\mathbf{e}_1 +\Omega =\Omega $, for any coherent $\mathcal{H}$-module $M$, and for any positive integer $i$, one has the following vanishing result:
\begin{equation}
\Ext^i_\mathcal{H}(M,C^\infty(\Omega ))=0.
\end{equation}
\end{thm}

%For complex domains, one might want to use Borel-Fourier transform instead of Fourier transform. In any case, one gets 

We can further prove the injectivity result for hyperfunctions by means of spectral sequence due to H. Komatsu and T. Oshima (cf. \cite{O}).

\begin{thm}
For any open convex subset $\Omega$ of $\mathbb{R}^n$ such that $\mathbb{R}\mathbf{e}_1 +\Omega =\Omega $, for any coherent $\mathcal{H}$-module $M$, and for any positive integer $i$, one has the following vanishing result:
\begin{equation}
\Ext^i_\mathcal{H}(M,\mathcal{B} (\Omega ))=0.
\end{equation}
\end{thm}

\begin{prf}
Put $U=\Omega\times\sqrt{-1}\mathbb{R}^n$, $U_i=\{ z\in U|\Image z_i\neq 0\}$, $\mathcal{U}=\{U,U_1\cdots ,U_n\}$, and $\mathcal{U}^\prime =\{U_1\cdots ,U_n\}$. By the Leray spectral sequence, one has 
\begin{equation}
H^p (\mathcal{U},\mathcal{U}^\prime ,\mathcal{O})=H^p_\Omega (U,\mathcal{O}).
\end{equation}
Now take a finite free resolution
\begin{equation}
0\rightarrow  \mathcal{H}^{1\times r_{N}}\xrightarrow{\times \mathbb{P}_N}\cdots \rightarrow \mathcal{H}^{1\times r_1}\xrightarrow{\times \mathbb{P}_1}\mathcal{H}^{1\times r_0}\rightarrow M\rightarrow 0\;(exact)
\end{equation}
of $M$. We put
\begin{align}
K^{p,q}&=C^p(\mathcal{U},\mathcal{U}^\prime ,\mathcal{O}^{r_q\times 1}),\\
\delta &=d^\prime :C^p(\mathcal{U},\mathcal{U}^\prime ,\mathcal{O}^{r_q\times 1})\rightarrow C^{p+1}(\mathcal{U},\mathcal{U}^\prime ,\mathcal{O}^{r_q\times 1}),\\
\mathbb{P}_q&=d^{\prime\prime }: C^p(\mathcal{U},\mathcal{U}^\prime ,\mathcal{O}^{r_q\times 1})\rightarrow C^p(\mathcal{U},\mathcal{U}^\prime ,\mathcal{O}^{r_{q+1}\times 1}).
\end{align}
By Theorem \ref{thm71} and the purity of relative cohomology (\cite{KKK} Chap. 2.), 
$$
\begin{aligned}
&^{\prime }E^{p,q}_1  =\;  ^{\prime\prime }H^q(K^{p,\cdot })  = \begin{cases}
                 C^p(\mathcal{U},\mathcal{U}^\prime ,\mathcal{O}^{r_0\times 1}_{\mathbb{P}_1}) & (p\leq n, q=0)\\
                 0 & (otherwise)
           \end{cases}
\\
&^{\prime\prime }E^{p,q}_1  =\;  ^{\prime }H^p(K^{\cdot ,q})  = \begin{cases}
                 \mathcal{B}(\Omega )^{r_q\times 1} & (p=n)\\
                 0 & (otherwise)
           \end{cases}
\\
&^{\prime }E^{p,q}_2  =\begin{cases}
                 H^p_\Omega (\mathcal{U},\mathcal{U}^\prime ,\mathcal{O}^{r_0\times 1}_{\mathbb{P}_1}) & (p\leq n, q=0)\\
                 0 & (otherwise)
           \end{cases}
\\
&^{\prime\prime }E^{p,q}_2  = \begin{cases}
\displaystyle
                 \frac{\Ker(\mathbb{P}_q:\mathcal{B}(\Omega )^{r_q\times 1}\rightarrow\mathcal{B}(\Omega )^{r_{q+1}\times 1})}{\Image(\mathbb{P}_{q-1}:\mathcal{B}(\Omega )^{r_{q-1}\times 1}\rightarrow\mathcal{B}(\Omega )^{r_{q}\times 1})} & (p=n)\\
                 0 & (otherwise).
           \end{cases}
\end{aligned}
$$
Here, $C^p(\mathcal{U},\mathcal{U}^\prime,\mathcal{O}_{\mathbb{P}_1}^{r_0\times 1})=\Ker\left( d^{\prime\prime }: C^p(\mathcal{U},\mathcal{U}^\prime ,\mathcal{O}^{r_0\times 1})\rightarrow C^p(\mathcal{U},\mathcal{U}^\prime ,\mathcal{O}^{r_{1}\times 1})\right).$ This shows that this spectral sequence degenerates at $E_2$ terms so that
\begin{equation}
0=\; ^{\prime }E^{n,q}_2=\; ^{\prime\prime }E^{n,q}_2=Ext^q_\mathcal{H}(M,\mathcal{B} (\Omega ))=0
\end{equation}
if $q>0.$
\qed
\end{prf}

We give another application of the division with bounds which is known as the problem of spectral synthesis (analysis).
We follow the argument of L. H\"ormander. The following lemma is taken from \cite{Ho}.

\begin{lem}[{\cite[Lemma 6.3.7]{Ho}}]\label{solvingeqs}
Let $\xi_1,\xi_2,\cdots$ be complex variables, $L_j\in\displaystyle\bigoplus_{i=1}^\infty\mathbb{C}\xi_i$ ($j=1,2,\cdots$), and $b=\{ b_j\}_{j\geq 1}\in\displaystyle\prod_{j=1}^\infty\mathbb{C}$. We put $\xi=\{\xi_i\}_i$. In these settings, an infinitely many linear equations $L_j(\xi )=b_j$  has a solution if and only if the following condition is satisfied:\\
\\
\;\;\;\;\;\;$\text{For any finite subset}\; F\subset\{1,2,\cdots\},\; \text{and for any complex numbers}$\\
\;\;\;\;\;\;$c_j\in\mathbb{C}\; (j\in F)\;\text{with}\; \displaystyle\sum_{j\in F}c_jL_j=0, \text{one has} \displaystyle\sum_{j\in F}c_jb_j=0$.
\end{lem}

\begin{thm}
Let $\mathscr{F}$ be either $C^\infty$ or $\mathcal{O}$ and  $\Omega$ be a convex subset of $\mathbb{R}^n$ (when $\mathscr{F}=C^\infty$) or of $\mathbb{C}^n$ (when $\mathscr{F}=\mathcal{O}$) such that $\Omega +\mathbb{R}\mathbf{e}_1$. For any  matrix of generalized $D\Delta$-operators $\mathbb{P}\in M(r_1,r_0;\mathcal{H})$, 
$\Ker\left(\mathbb{P}(\partial ,\sigma ):\mathbf{E}^{r_0\times 1}\rightarrow\mathbf{E}^{r_1\times 1}\right)$ is dense in $\Ker(\mathbb{P}(\partial ,\sigma ):\mathscr{F}(\Omega)^{r_0\times 1}\rightarrow\mathscr{F}(\Omega)^{r_1\times 1}).$
\end{thm}

\begin{prf}
%Since the topology of $\mathcal{O}(\Omega )$ is induced from $C^\infty(\Omega )$, we can assume $\mathscr{F}=C^\infty$.
We prove the theorem only for $\mathscr{F}=C^\infty$ since the arguments are similar.
First, note that we have the identity
\begin{equation}\overline{\Ker(\mathbb{P}(\partial ,\sigma ):\mathbf{E}^{r_0\times 1}\rightarrow \mathbf{E}^{r_1\times 1})}=\bigcap\{\Ker L\mid L\in (C^\infty(\Omega)^{r_0\times 1})^\prime ,\;L\restriction_V=0\}
\end{equation}
by the Hahn-Banach theorem, where the bar stands for the closure and  
\begin{equation}
V=\Ker(\mathbb{P}(\partial ,\sigma ):\mathbf{E}^{r_0\times 1}\rightarrow\mathbf{E}^{r_1\times 1}).
\end{equation}
Let $L$ be a linear functional $L= \;{}^t(f_1,\cdots ,f_{r_0})\in (C^\infty(\Omega)^{r_0\times 1})^\prime$ with $L|_V=0$. For this $L$, there exists a compact convex subset $K$ of $\Omega$, a positive real number $C>0,$ and  a non-negative integer $M$  such that, for all $\zeta\in\mathbb{C}^n,$
\begin{equation}
|\widehat{L}(\zeta )|\leq C(1+|\zeta|)^Me^{H_K(\zeta)},
\end{equation}
where $\widehat{L}$ is the Fourier transform of $L$.

Now suppose we could find an element $\mathbf{g}\in (C^\infty(\Omega)^{r_1\times 1})^\prime $ such that \\
$L= ^{t}\mathbb{P}(-\partial ,\sigma^{-1})\mathbf{g}$. 
Then, for any $\mathbf{u}\in C^\infty(\Omega)^{r_0\times 1}$,
\begin{equation}
\langle L,\mathbf{u}\rangle=\langle ^{t}\mathbb{P}(-\partial ,\sigma^{-1})\mathbf{g} ,\mathbf{u}\rangle=\langle\mathbf{g},\mathbb{P}(\partial,\sigma)\mathbf{u}\rangle.
\end{equation}
This implies $\langle L,\mathbf{u}\rangle=0$ if $\mathbb{P}(\partial,\sigma)\mathbf{u}=0$ in a neighbourhood of $\supp \mathbf{g}$, hence we can get the theorem.

Now the proof of this theorem is reduced to solving the equation
\begin{equation}\label{FundamentalEquation}
\widehat{L}= \mathbb{Q}(\zeta)\hat{\mathbf{g}}(\zeta)
\end{equation}
with $\hat{\mathbf{g}}\in\mathcal{O}(\mathbb{C}^n)_p^{r_1\times 1}$ via the Fourier transform, where $\mathbb{Q}(\zeta)$ is a matrix of holomorphic functions defined by 
$\mathbb{Q}(\zeta)= ^{t}\mathbb{P}(-\sqrt{-1}\zeta,e^{-\sqrt{-1}\zeta_1})$. The proof falls naturally into two parts.\\
$\underline{Step1}$ Construction of the formal solution.

For any fixed $\zeta_0\in\mathbb{C}^n$, we are going to construct a formal series solution $\hat{\mathbf{g}}=G_{\zeta_0}\in\mathbb{C}[[\zeta-\zeta_0]]$ of (\ref{FundamentalEquation}). Note that the equation (\ref{FundamentalEquation}) in this case can be written as an infinite system of equations
\begin{equation}
\partial^\alpha_\zeta\widehat{L}(\zeta)\restriction_{\zeta=\zeta_0}=\partial^\alpha_\zeta\left(\mathbb{Q}(\zeta)G_{\zeta_0}(\zeta)\right)\restriction_{\zeta=\zeta_0},\;\;(\alpha\in\mathbb{Z}^n_{\geq 0}).
\end{equation}

\noindent
By Lemma \ref{solvingeqs}, it is enough to show the following: Take a polynomial vector $\mathbf{q}\in\mathbb{C}[\partial]^{1\times r_0}$. If the equation
\begin{equation}\label{CompatibilityEq}
\mathbf{q}(\partial_\zeta)\cdot\left(\mathbb{Q}(\zeta)\mathbf{h}(\zeta)\right)\restriction_{\zeta=\zeta_0}=0
\end{equation}
holds for all formal series $\mathbf{h}(\zeta)\in\mathbb{C}[[\zeta-\zeta_0]]^{r_1\times 1}$, then one has
\begin{equation}
\left(\mathbf{q}(\partial_\zeta)\cdot\widehat{L}(\zeta)\right)\restriction_{\zeta=\zeta_0}=0.
\end{equation}
Take any $\mathbf{q}=(q_1,\cdots ,q_{r_0})\in\mathbb{C}[\partial]^{1\times r_0}$ which satisfies the equation (\ref{CompatibilityEq}) and put $u_k(z)=q_k(-\sqrt{-1}z)e^{-\sqrt{-1}z\zeta_0},\;\mathbf{u}= {}^{t}(u_1(z),\cdots,u_{r_0}(z))$.
We have 
\begin{align}
\mathbb{P}(\partial_z,\sigma_z)\cdot\mathbf{u}&=\mathbb{P}(\partial_z,\sigma_z)\cdot \left({}^{t}(q_1(\partial_\zeta),\cdots,q_{r_0}(\partial_\zeta))\cdot e^{-\sqrt{-1}z\zeta}\right)\restriction_{\zeta=\zeta_0}\\
&={}^t\left\{(q_1(\partial_\zeta),\cdots,q_{r_0}(\partial_\zeta))\cdot\left({}^t\mathbb{P}(-\sqrt{-1}\zeta,e^{-\sqrt{-1}\zeta_1})e^{-\sqrt{-1}z\zeta}\right)\right\}\restriction_{\zeta=\zeta_0}\\
&=0.
\end{align}
Thus, $\mathbf{u}(z)$ is an exponential polynomial solution. Since $L\restriction_V=0$, we have 
\begin{equation}
0=\langle L,\mathbf{u}\rangle=\left(\langle L,(q_1(\partial_\zeta),\cdots,q_{r_0}(\partial_\zeta))e^{-\sqrt{-1}z\zeta})\rangle_{z}\right)\restriction_{\zeta=\zeta_0}=\mathbf{q}(\partial_\zeta)\cdot\widehat{L}(\zeta)\restriction_{\zeta=\zeta_0},
\end{equation}
and hence we get the formal solution $G_{\zeta_0}$.\\
$\underline{Step2}$ Construction of the solution with bounds

Since $\Omega$ is Stein, global section functor $\Gamma(\Omega ,-)$ is exact for complexes of coherent sheaves. We get from Step1 that $\widehat{L}\in\mathbb{Q}(\zeta)\cdot\mathcal{O}(\mathbb{C}^n)^{r_1\times 1}$. Now we can apply Theorem \ref{division} to $\widehat{L}$ to obtain the desired solution $\hat{\mathbf{g}}\in\mathcal{O}(\mathbb{C}^n)^{r_1\times 1}_p$.\qed
\end{prf}

Lastly, we give a description of cohomology groups of the solution sheaves. 
Let $\pi:\mathbb{R}^n\rightarrow \mathbb{R}^{n-1}$ be the projection which truncates the first coordinate. For any given matrix of $D\Delta$-operators $\mathbb{P}\in M(r_1,r_0;\mathcal{H}),$ we can define the subsheaf $\mathcal{S}ol_{\mathbb{P},\mathscr{F}}$ of $\pi_*\mathscr{F}^{r_0\times 1}$ for $\mathscr{F}=C^\infty ,\mathcal{B}$ by the formula
\begin{equation}
\mathcal{S}ol_{\mathbb{P},\mathscr{F}}= \Ker\left(\mathbb{P}:\pi_*\mathscr{F}^{r_0\times 1}\rightarrow\pi_*\mathscr{F}^{r_1\times 1}\right).
\end{equation}

\begin{thm}
We put $M=\mathcal{H}^{1\times r_0}/\mathcal{H}^{1\times r_1}\cdot\mathbb{P}.$ For any open set $\Omega^\prime\subset\mathbb{R}^{n-1}$, there is a canonical cohomology isomorphism
\begin{equation}
H^{i}(\Omega^\prime ,\mathcal{S}ol_{\mathbb{P} ,\mathscr{F}})\simeq \Ext^i_{\mathcal{H}}(M, \pi_*\mathscr{F}(\Omega^\prime )),
\end{equation}
where $\mathscr{F}=C^\infty, \mathcal{B}$ and $i$ is any integer.

\end{thm}

\begin{prf}

We first prove the theorem for $\mathscr{F}=C^\infty.$ Consider a finite free resolution 
\begin{equation}\label{Formula638}
0\rightarrow  \mathcal{H}^{1\times r_{N}}\xrightarrow{\times \mathbb{P}_N}\cdots \rightarrow \mathcal{H}^{1\times r_1}\xrightarrow{\times \mathbb{P}_1=\times\mathbb{P}}\mathcal{H}^{1\times r_0}\rightarrow M\rightarrow 0\;(exact)
\end{equation}
of $M.$ Applying the functor $\mathcal{H}om_{\mathcal{H}}(-,\pi_*\mathscr{F})$ to this sequence, we have a complex
\begin{equation}
0\rightarrow \mathcal{S}ol_{\mathbb{P} ,\mathscr{F}}\rightarrow\pi_*\mathscr{F}^{r_0\times 1}\xrightarrow{\mathbb{P}_1\cdot } \pi_*\mathscr{F}^{r_1\times 1}\cdots \xrightarrow{\mathbb{P}_N\cdot } \pi_*\mathscr{F}^{r_N\times 1}\rightarrow 0.
\end{equation}

\noindent
In view of Theorem \ref{thm72} and the fact that open convex subsets are fundamental system of neighbourhoods in $\mathbb{R}^{n-1},$ we can conclude that the complex above is exact. Furthermore, since for any $i>0,$ we have
\begin{equation}
H^{i}(\Omega^\prime,\pi_*\mathscr{F})=0,
\end{equation}
the exact sequence above is a $\Gamma (\Omega^\prime,-)$-injective resolution of $\mathcal{S}ol_{\mathbb{P},\mathscr{F}}.$
Therefore, by the standard theory of sheaf cohomology, we obtain
\begin{equation}
H^i(\Omega^\prime,\mathcal{S}ol_{\mathbb{P},\mathscr{F}})\simeq\frac{\Ker\left(\mathbb{P}_i\cdot :\pi_*\mathscr{F}(\Omega^\prime)^{r_i\times 1}\rightarrow \pi_*\mathscr{F}(\Omega^\prime)^{r_{i+1}\times 1} \right)}{\Image\left(\mathbb{P}_{i-1}\cdot :\pi_*\mathscr{F}(\Omega^\prime)^{r_{i-1}\times 1}\rightarrow \pi_*\mathscr{F}(\Omega^\prime)^{r_{i}\times 1} \right)}.
\end{equation}
On the other hand, if one applies the functor $\Hom_\mathcal{H}(-,\pi_*\mathscr{F}(\Omega^\prime))$ to the finite free resolution (\ref{Formula638}) of $M,$ one has
\begin{equation}
\Ext^i_{\mathcal{H}}(M, \pi_*\mathscr{F}(\Omega^\prime ))\simeq\frac{\Ker\left(\mathbb{P}_i\cdot :\pi_*\mathscr{F}(\Omega^\prime)^{r_i\times 1}\rightarrow \pi_*\mathscr{F}(\Omega^\prime)^{r_{i+1}\times 1} \right)}{\Image\left(\mathbb{P}_{i-1}\cdot :\pi_*\mathscr{F}(\Omega^\prime)^{r_{i-1}\times 1}\rightarrow \pi_*\mathscr{F}(\Omega^\prime)^{r_{i}\times 1} \right)}.
\end{equation}

As for $\mathscr{F}=\mathcal{B}$, the argument is similar since 
\begin{equation}
H^{i}(\Omega^\prime,\pi_*\mathscr{F})=0
\end{equation}
for $i>0$ follows from the flabbiness of $\mathcal{B}.$
\qed
\end{prf}




\section*{Acknowledgement}

The author would like to thank H. Sakai and T. Oshima for their constant encouragement and uncountably many suggestions, to Y. Nozaki and S. Wakatsuki for various comments, and to anonymous referees whose comments improved the exposition of this paper.

\end{document}